\documentclass[a4paper]{amsart}

\theoremstyle{plain}

\newtheorem{theorem}{Theorem}[section]
\newtheorem{proposition}[theorem]{Proposition}
\newtheorem{lemma}[theorem]{Lemma}
\newtheorem{corollary}[theorem]{Corollary}

\theoremstyle{definition}

\newtheorem{definition}[theorem]{Definition}

\newtheorem{remark}[theorem]{Remark}

\RequirePackage{amssymb}

\RequirePackage[all]{xy}

\UseComputerModernTips    

\entrymodifiers={!!<0ex,.7ex>+}

\newcommand{\isoar}{\ar|-*=0!/_1.5pt/\dir{~}}

\newcommand{\ooint}[1]{\left]{#1}\right[}

\newcommand{\union}{\cup}
\newcommand{\Union}{\bigcup}
\newcommand{\inter}{\cap}

\newcommand{\N}{\mathbb{N}}
\newcommand{\Z}{\mathbb{Z}}

\newcommand{\C}{\mathbb{C}}

\newcommand{\scal}[1]{\left\langle{#1}\right\rangle}
\newcommand{\norm}[1]{\left\Vert{#1}\right\Vert}

\newcommand{\op}{\mathrm{op}}

\DeclareMathOperator{\id}{id}

\newcommand{\DSum}{\bigoplus}
\newcommand{\dsum}[1][]{\mathbin{\oplus_{#1}}}

\newcommand{\indlim}{\mathop{\varinjlim}\limits}

\newcommand{\prolim}{\mathop{\varprojlim}\limits}

\DeclareMathOperator{\im}{im}
\DeclareMathOperator{\coim}{coim}

\newcommand{\comp}{\mathbin{\circ}}

\renewcommand{\to}[1][]{\xrightarrow{#1}}

\newcommand{\isoto}{\xrightarrow{\sim}}

\newcommand{\Hom}[1][]{\mathrm{Hom}_{\raise1.5ex\hbox to.1em{}#1}}

\newcommand{\RHom}[1][]{\mathrm{RHom}_{\raise1.5ex\hbox to.1em{}#1}}

\def\shf{\mathcal{F}}
\def\shg{\mathcal{G}}

\def\shp{\mathcal{P}}

\def\shr{\mathcal{R}}

\renewcommand{\hom}[1][]{{\mathcal{H}om}_{\raise1.5ex\hbox to.1em{}#1}}

\newcommand{\rhom}[1][]{{R\mathcal{H}om}_{\raise1.5ex\hbox to.1em{}#1}}

\newcommand{\ext}[2][]{%
{\mathcal{E}xt}_{\raise1.5ex\hbox to.1em{}#1}^{#2}}

\newcommand{\tens}[1][]{%
\mathbin{\otimes_{\raise1.5ex\hbox to-.1em{}#1}}}

\newcommand{\ltens}[1][]{%
\mathbin{\otimes_{\raise1.5ex\hbox to-.1em{}#1}^{L}}}

\newcommand{\eim}[1]{{#1}_!}

\newcommand{\reim}[1]{{R#1}_!}

\newcommand{\opb}[1]{#1^{-1}}

\newcommand{\epb}[1]{#1^{!}}

\newcommand{\tenstop}[1][]{%
\mathbin{\hat{\otimes}_{\raise1.5ex\hbox to-.1em{}#1}}}

\newcommand{\homtop}[1][]{\mathcal{L}_{\raise1.5ex\hbox to.1em{}#1}}

\newcommand{\Homtop}[1][]{\mathrm{L}_{\raise1.5ex\hbox to.1em{}#1}}

\newcommand{\D}{\mathcal{D}}

\renewcommand{\O}{\mathcal{O}}

\newcommand{\dsol}{\mathop{\mathcal{S}ol}\nolimits}

\newcommand{\lh}{\mathcal{LH}}

\newcommand{\ctb}{\mathcal{B}}
\newcommand{\ctc}{\mathcal{C}}

\newcommand{\cte}{\mathcal{E}}

\newcommand{\cti}{\mathcal{I}} 
\newcommand{\ctj}{\mathcal{J}}

\newcommand{\ctm}{\mathcal{M}}
\newcommand{\ctn}{\mathcal{N}} 
\newcommand{\ctr}{\mathcal{R}}

\newcommand{\ctu}{\mathcal{U}}
\newcommand{\ctv}{\mathcal{V}}

\newcommand{\ctBan}{\mathcal{B}an}
\newcommand{\ctFr}{\mathcal{F}r}
\newcommand{\ctTc}{\mathcal{T}c}
\newcommand{\ctAb}{\mathcal{A}b} 
\newcommand{\Ind}{\mathcal{I}nd}
\newcommand{\Pro}{\mathcal{P}ro}

\newcommand{\rprolim}{\mathrm{R}\prolim}
\newcommand{\lindlim}{\mathrm{L}\indlim}
\newcommand{\cproobj}[1]{\text{``}{#1}\text{''}}
\newcommand{\cindobj}[1]{\text{``}{#1}\text{''}}
\newcommand{\fprolim}{%
\mathop{\text{``}\mathord{\varprojlim}\text{''}}\limits}
\newcommand{\findlim}{%
\mathop{\text{``}\mathord{\varinjlim}\text{''}}\limits}

\DeclareMathOperator{\iolim}{IL}

\DeclareMathOperator{\IB}{IB}
\DeclareMathOperator{\Sect}{\Gamma}
\DeclareMathOperator{\RSect}{R\Gamma}

\newcommand{\Dual}{\mathrm{D}}
\newcommand{\IDual}{\mathrm{D}_{\mathrm{i}}}
\newcommand{\SDual}{\mathrm{D}_{\mathrm{b}}}

\newcommand{\ltenstop}[1][]{\mathbin{\hat{\otimes}_{\raise1.5ex\hbox 
to-.1em{}#1}^{L}}}
\newcommand{\RHomtop}[1][]{\mathrm{RL}_{\raise1.5ex\hbox to.1em{}#1}}
\newcommand{\rhomtop}[1][]{%
\mathrm{R}\mathcal{L}_{\raise1.5ex\hbox to.1em{}#1}}
\newcommand{\etenstop}{\hat{\boxtimes}}
\newcommand{\letenstop}{\hat{\boxtimes}^{L}}

\newcommand{\Shv}{%
\mathop{\mathcal{S}\mspace{-2mu}\mathit{hv}}\nolimits}
\newcommand{\ctMod}{\mathcal{M}od}

\newcommand{\Uu}{\mathfrak{U}}
\newcommand{\Uv}{\mathfrak{V}}

\begin{document}

\title{A Topological Reconstruction Theorem for $\D^{\infty}$-Modules}

\author{Fabienne Prosmans}

\address{%
Laboratoire Analyse, G\'{e}om\'{e}trie et Applications\\
UMR 7539\\
Universit\'{e} Paris 13\\
Avenue J.-B. Cl\'{e}ment\\
93430 Villetaneuse\\
France}

\email{prosmans@math.univ-paris13.fr}

\author{Jean-Pierre Schneiders}

\address{%
Laboratoire Analyse, G\'{e}om\'{e}trie et Applications\\
UMR 7539\\
Universit\'{e} Paris 13\\
Avenue J.-B. Cl\'{e}ment\\
93430 Villetaneuse\\
France}

\email{jps@math.univ-paris13.fr}

\date{}



\begin{abstract}
In this paper, we prove that any perfect complex of 
$\D^{\infty}$-modules may be reconstructed from its holomorphic 
solution complex provided that we keep track of the natural topology 
of this last complex.  This is to be compared with the reconstruction 
theorem for regular holonomic $\D$-modules which follows from the 
well-known Riemann-Hilbert correspondence.  To obtain our result, we 
consider sheaves of holomorphic functions as sheaves with values in 
the category of ind-Banach spaces and study some of their homological 
properties.  In particular, we prove that a K\"{u}nneth formula holds 
for them and we compute their Poincar\'{e}-Verdier duals.  As a 
corollary, we obtain the form of the kernels of ``continuous'' 
cohomological correspondences between sheaves of holomorphic forms.  
This allows us to prove a kind of holomorphic Schwartz' kernel theorem 
and to show that 
$\D^{\infty}\simeq\mathrm{R}\mathcal{H}\mathit{omtop}(\O,\O)$.  Our 
reconstruction theorem is a direct consequence of this last 
isomorphism.  Note that the main problem is the vanishing of the 
topological Ext's and that this vanishing is a consequence of the 
acyclicity theorems for DFN spaces that are established in the paper.
\end{abstract}

\maketitle

\tableofcontents

\setcounter{section}{-1}

\section{Introduction}

In algebraic analysis, one represents systems of analytic linear 
partial differential equations on a complex analytic manifold $X$ by 
modules over the ring $\D_{X}$ of linear partial differential 
operators with analytic coefficients.  Using this representation, the 
holomorphic solutions of the homogeneous system associated to the 
$\D_{X}$-module $\ctm$ correspond to
$$
\hom[\D_{X}](\ctm,\O_{X})
$$
where $\O_{X}$ denotes the $\D_{X}$-module of holomorphic functions. 
If one wants also to take into consideration the compatibility 
conditions, one has to study the full solution complex
$$
\dsol(\ctm)=\rhom[\D_{X}](\ctm,\O_{X})
$$
in the derived category $D^{+}(\C_{X})$ of sheaves of $\C$-vector 
spaces.  In~\cite{Kashiwara84b} (see also~\cite{Mebkhout84b}), it was 
shown that the functor $\dsol$ induces an equivalence between the 
derived category formed by the bounded complexes of regular holonomic 
$\D_{X}$-modules and that formed by the bounded complexes of 
$\C$-constructible $\C_{X}$-modules.  This equivalence is usually 
called the Riemann-Hilbert correspondence.  One of its corollaries is 
that it is possible to reconstruct a complex of regular holonomic 
$\D_{X}$-modules from its complex of holomorphic solutions.

Our aim in this paper is to extend this reconstruction theorem to 
perfect complexes of $\D_{X}^{\infty}$-modules by taking into account 
the natural topology of the complex of holomorphic solutions. 
Informally, the relation we will obtain is of the type
$$
\ctm\simeq\mathrm{R}\mathcal{H}\mathit{omtop}(\dsol(\ctm),\O_{X})
$$
and will follow from the fact that 
$$
\D_{X}^{\infty}\simeq
\mathrm{R}\mathcal{H}\mathit{omtop}(\O_{X},\O_{X}).
$$
To give a meaning to these formulas, we will have to work in the 
derived category of sheaves with values in the category of ind-objects 
of the category of Banach spaces using the techniques and results 
of~\cite{Schneiders99}.

As a help to the reader, let us briefly recall the main facts 
concerning quasi-abelian homological algebra and sheaf theory 
established in that paper.

The central notion is that of a quasi-abelian category (i.e.\ an 
additive category with kernels and cokernels such that the 
push-forward (resp.\ the pull-back) of a kernel (resp.\ a cokernel) is 
still a kernel (resp.\ a cokernel)). Let $\cte$ be such a category. A 
morphism of $\cte$ is said to be strict if its coimage is canonically 
isomorphic to its image and a complex
$$
\cdots \to X^{k-1} \to[d^{k-1}] X^{k} \to[d^{k}] X^{k+1} \to \cdots
$$
of $\cte$ is said to be strictly exact in degree $k$ if $d^{k-1}$ is 
strict and $\ker d^{k} = \im d^{k-1}$.  Localizing the triangulated 
category $K(\cte)$ of complexes ``modulo homotopy'' by the null system 
formed by the complexes which are strictly exact in every degree gives 
us the derived category $D(\cte)$.  This category has two canonical 
t-structures.  Here, we will only use the left one.  Its heart 
$\lh(\cte)$ is formed by the complexes of the form
$$
0 \to X^{-1} \to[d^{-1}] X^{0} \to 0
$$
where $d^{-1}$ is a monomorphism. The cohomology functor
$$
LH^{k} : D(\cte) \to \lh(\cte)
$$
sends the complex $X^{\cdot}$ to the complex
$$
0\to \coim d^{k-1} \to \ker d^{k} \to 0
$$
with $\ker d^{k}$ in degree $0$.  In~\cite{Schneiders99}, it was shown 
that in most problems of homological algebra and sheaf theory we may 
replace the quasi-abelian category $\cte$ by the abelian category 
$\lh(\cte)$ without loosing any information.  It was also shown there 
that if $\cte$ is elementary (i.e.\ if it has a small strictly 
generating set formed by tiny projective objects), then the sheaves 
with values in $\cte$ share most of the usual properties of sheaves of 
abelian groups (including Poincar\'{e}-Verdier duality).  If $\cte$ 
has moreover a closed structure given by an internal tensor product 
and an internal $\Hom$ functor satisfying some natural assumptions 
then K\"{u}nneth theorem holds for sheaves with values in $\cte$.

Let us now introduce the quasi-abelian categories that will be used 
in this paper and fix our notations.

Following~\cite{Prosmans99a}, we denote $\ctBan$ (resp.\ $\ctFr$, 
$\ctTc$) the quasi-abelian category of Banach spaces (resp.\ 
Fr\'{e}chet spaces, arbitrary locally convex topological vector 
spaces).  Let us recall (see e.g.~\cite{Prosmans95}) that, for any set 
$I$, the space $l^{1}(I)$ (resp.\ $l^{\infty}(I)$) of summable (resp.\ 
bounded) sequences of $\C$ indexed by $I$ is projective (resp.\ 
injective) in $\ctBan$.  Using these spaces, one shows easily that 
$\ctBan$ has enough injective and projective objects.  Recall also 
that the category $\ctBan$ has a canonical structure of closed 
additive category given by a right exact tensor product
$$
\tenstop:\ctBan\times\ctBan\to\ctBan
$$
and a left exact internal $\Hom$
$$
\Homtop:\ctBan^{\op}\times\ctBan\to\ctBan.
$$
Denoting $\ltenstop$ the left derived functor of $\tenstop$ and 
$\RHomtop$ the right derived functor of $\Homtop$, we have the 
adjunction formula
$$
\RHom(E\ltenstop F,G)\simeq
\RHom(E,\RHomtop(F,G)).
$$

Let $\Uu$, $\Uv$ be two universes such that $\Uv\ni\Uu$. As usual, 
denote $\ctBan_{\Uu}$ the category formed by the Banach spaces which 
belong to $\Uu$ and consider the category
$$
\Ind_{\Uv}(\ctBan_{\Uu})
$$
of ind-objects of $\ctBan_{\Uu}$. Recall that the objects of 
$\Ind_{\Uv}(\ctBan_{\Uu})$ are functors
$$
E:\cti\to\ctBan_{\Uu}
$$
where $\cti$ is a $\Uv$-small filtering category and that if
$$
E:\cti\to\ctBan_{\Uu},\qquad
F:\ctj\to\ctBan_{\Uu}
$$
are two such functors, then
$$
\Hom[\Ind_{\Uv}(\ctBan_{\Uu})](E,F)=
\prolim_{i\in\cti}\indlim_{j\in\ctj}\Hom[\ctBan_{\Uu}](E(i),F(j)).
$$
For further details on ind-objects, we refer the reader to classical 
sources (such as~\cite{SGA4a,Artin-Mazur69}) and 
to~\cite{Prosmans98a}.  Following the standard usage and to avoid 
confusions, we will denote
$$
\findlim_{i\in\cti}E(i)
$$
the functor $E:\cti\to\ctBan_{\Uu}$ 
considered as an object of $\Ind_{\Uv}(\ctBan_{\Uu})$. Similarly, we 
denote $\cindobj{X}$ the ind-object associated to the $\Uu$-Banach 
space $X$. In other words, we set
$$
\cindobj{X}=\findlim_{i\in\cti}C(i)
$$
where $\cti$ is a one point category and $C:\cti\to\ctBan_{\Uu}$ is 
the constant functor with value $X$.  Note that, in the rest of the 
paper, we will not make the universes $\Uu$, $\Uv$ explicit in our 
notations since this is not really necessary for a clear 
understanding.

Using~\cite{Schneiders99}, we see that the category $\Ind(\ctBan)$ of 
ind-objects of $\ctBan$ is an elementary closed quasi-abelian 
category.  It follows that sheaves with values in $\Ind(\ctBan)$ share 
most of the usual properties of abelian sheaves (including K\"{u}nneth 
Theorem and Poincar\'{e}-Verdier duality).  In $\Ind(\ctBan)$, the 
internal tensor product
$$
\tenstop:\Ind(\ctBan)\times\Ind(\ctBan)\to\Ind(\ctBan)
$$
and the internal $\Hom$ functor
$$
\Homtop:(\Ind(\ctBan))^{\op}\times\Ind(\ctBan)\to\Ind(\ctBan)
$$
are characterized by
$$
(\findlim_{i\in I}E_{i})\tenstop(\findlim_{j\in J}F_{j})=
\indlim_{i\in I}\indlim_{j\in J}\cindobj{E_{i}\tenstop F_{j}}
$$
and
$$
\Homtop(\findlim_{i\in I}E_{i},\findlim_{j\in J}F_{j})=
\prolim_{i\in I}\indlim_{j\in J}\cindobj{\Homtop(E_{i},F_{j})}.
$$
The internal tensor product (resp.\ internal $\Hom$ 
functor, external tensor product) for sheaves with values in 
$\Ind(\ctBan)$ will be denoted by $\tenstop$ (resp.\ $\homtop$, 
$\etenstop$) and we will use the other notations of sheaf theory in 
their usual form. In particular, $\reim{f}$ (resp.\ $\epb{f}$) will 
denote the direct (resp.\ inverse) image with proper support and 
$\omega_{X}$ (resp.\ $\Dual(\cdot)=\rhomtop(\cdot,\omega_{X})$) will 
denote the Poincar\'{e}-Verdier dualizing complex (resp.\ functor).

Let us now describe with some details the content of this paper.

In the first section, we study the functor $\IB:\ctTc\to\Ind(\ctBan)$ 
defined by setting
$$
\IB(E)=\findlim_{B\in\ctb_{E}}\widehat{E}_{B}
$$
where $\ctb_{E}$ is the set of absolutely convex bounded subsets of 
$E$ and $E_{B}$ the linear hull of $B$.  We establish the properties 
of this functor we need in the rest of the paper.  More precisely, we 
prove that if $E$ is bornological and $F$ complete, then
$$
\Hom[\Ind(\ctBan)](\IB(E),\IB(F))\simeq
\Hom[\ctTc](E,F)
$$
and
$$
\IB(\Homtop[b](E,F))\simeq\Homtop(\IB(E),\IB(F)).
$$
Here, $\Homtop[b](E,F)$ is the vector space $\Hom[\ctTc](E,F)$ endowed 
with the system of semi-norms
$$
\{p_{B}:p\text{ continuous semi-norm of }F, 
B\text{ bounded subset of }E\}
$$
where
$$
p_{B}(h)=\sup_{e\in B}p(h(e)).
$$
Moreover, we show that $\IB$ is compatible with projective limits of 
filtering projective systems of complete spaces.  We show also its 
compatibility with complete inductive limits of injective inductive 
systems of Fr\'{e}chet spaces indexed by $\N$.

The second section is devoted to the proof of some acyclicity results 
for $\Homtop$ and $\tenstop$ in $\Ind(\ctBan)$.  First, we show that 
if $E$ is a DFN space and if $F$ is a Fr\'{e}chet space, then both 
$LH^{k}(\RHom(\IB(E),\IB(F)))$ and $LH^{k}(\RHomtop(\IB(E),\IB(F)))$ 
are $0$ for $k\neq 0$.  (Note that a related result was obtained for 
the category $\ctTc$ by Palamodov in~\cite{Palamodov71}.)  Next, we 
establish that if $E$ and $F$ are objects of $\Ind(\ctBan)$ with $E$ 
nuclear, then
\begin{equation}\tag{*}
E\ltenstop F\simeq E\tenstop F.
\end{equation}

We start Section~3 by proving that if $X$ is a topological 
space with a countable basis and if $F$ is a presheaf of Fr\'{e}chet 
spaces on $X$ which is a sheaf of vector spaces, then
$$
U\mapsto\IB(F(U))
\qquad(U\text{ open of }X)
$$
is a sheaf with values in $\Ind(\ctBan)$.  This shows, in particular,
that $\IB(\O_{X})$ is a sheaf with values in $\Ind(\ctBan)$ for any 
complex analytic manifold $X$.  We end the section by establishing 
that
$$
\RSect(U,\IB(\O_{X}))\simeq
\Sect(U,\IB(\O_{X}))
$$
if $U$ is an open subset of $X$ such that $H^{k}(U,\O_{X})\simeq 0$ 
($k\neq 0$).  This result may be viewed as a topological version of 
Cartan's Theorem~B. As a corollary, if $X$ is a Stein manifold, we get 
a similar isomorphism with $U$ replaced by any holomorphically convex 
compact subset of $X$.

In Section~4, using (*), we show that
$$
\IB(\O_{X})\letenstop\IB(\O_{Y})
\simeq\IB(\O_{X\times Y})
$$
for any complex analytic manifolds $X$ and $Y$. This allows us to 
obtain a topological K\"{u}nneth Theorem for holomorphic cohomology.

Section~5 is devoted to the proof that, for any complex analytic 
manifold $X$ of dimension $d_{X}$, the Poincar\'{e} dual of 
$\IB(\O_{X})$ is isomorphic to $\IB(\Omega_{X})[d_{X}]$.  Since the 
problem is of local nature, we find, by a series of reductions using 
the results established in the previous sections, that it is 
sufficient to show that, if $P$ is a closed interval of $\C$ and $V$ 
is an open interval of $\C^{n}$, then
$$
\RSect_{P\times V}(\C\times V,\IB(\O_{\C\times V}))
\simeq\Homtop(\IB(\O_{\C}(P)),\IB(\O_{V}(V)))[-1].
$$
This isomorphism is obtained by proving that, in this situation, one 
has a split exact sequence of the form
$$
0\to\O_{\C\times V}(\C\times V)\to
\O_{\C\times V}((\C\setminus P)\times V)\to
\Homtop[b](\O_{\C}(P),\O_{V}(V))\to 0
$$
in $\ctTc$.

We begin Section~6 by giving the general form of the kernels of 
continuous cohomological correspondences between sheaves of 
holomorphic differential forms.  More precisely, we show that, if $X$, 
$Y$ are complex analytic manifolds of dimension $d_{X}$, $d_{Y}$, then
$$
\IB(\Omega_{X\times Y}^{(d_{X}-r,s)})[d_{X}]\simeq
\rhomtop(q_{X}^{-1}\IB(\Omega_{X}^{r}),\epb{q_{Y}}\IB(\Omega_{Y}^{s})).
$$
As a consequence, we find that, for any morphism of complex analytic 
manifolds $f:X\to Y$, we have a canonical isomorphism
$$
\rhomtop(\opb{f}\IB(\O_{Y}),\IB(\O_{X}))\simeq
\opb{\delta_{f}}\RSect_{\Delta_{f}}
\IB(\Omega_{X\times Y}^{(0,d_{Y})})[d_{Y}]
$$
where $\Delta_{f}$ is the graph of $f$ in $X\times Y$ and 
$\delta_{f}:X\to X\times Y$ is the associated graph embedding.  In 
particular,
$$
LH^{k}(\rhomtop(\opb{f}\IB(\O_{Y}),\IB(\O_{X})))=0
$$
for $k\neq 0$ and
\begin{equation}\tag{**}
\rhom(\opb{f}\IB(\O_{Y}),\IB(\O_{X}))\simeq
\D_{X\to Y}^{\infty}.
\end{equation}
Note that this contains the fact that continuous endomorphisms of 
$\O_{X}$ may be identified with partial differential operators of 
infinite order as was conjectured by Sato and established by Ishimura 
in~\cite{Ishimura78}.

We start the last section by proving an abstract reconstruction 
theorem for perfect complexes of modules over a ring in the closed 
category $\Shv(X;\Ind(\ctBan))$. Thanks to the embedding functor
$$
\tilde{I}_{\ctv}:\Shv(X;\ctv)\to\Shv(X;\Ind(\ctBan))
$$
(where $\ctv$ denotes the category of $\C$-vector spaces) we are also 
able to prove a similar formula for perfect complexes of modules over 
an ordinary sheaf of rings.  Using (**) with $f=\id_{X}$, we get a 
topological reconstruction theorem for $\D_{X}^{\infty}$-modules.  
More precisely, we prove that the functors
$$
\rhomtop[\tilde{I}_{\ctv}(\D_{X}^{\infty})]
(\tilde{I}_{\ctv}(\cdot),\IB(\O_{X})):
D^{-}(\ctMod(\D_{X}^{\infty}))
\to D^{+}(\Shv(X;\Ind(\ctBan)))
$$
and 
$$
\rhom(\cdot,\IB(\O_{X})):
D^{-}(\Shv(X;\Ind(\ctBan)))\to
D^{+}(\ctMod(\D_{X}^{\infty})) 
$$
are well-defined and that
$$
\rhom(\rhomtop[\tilde{I}_{\ctv}(\D_{X}^{\infty})]
(\tilde{I}_{\ctv}(\ctm),\IB(\O_{X})),\IB(\O_{X}))\simeq\ctm
$$
for any perfect complex of $\D_{X}^{\infty}$-modules $\ctm$.  Note 
that the image of $\ctm$ by the first functor above is a kind of 
topologized version of the holomorphic solution complex of $\ctm$ and 
that the preceding formula may be viewed as a way to reconstruct a 
perfect system of analytic partial differential equations of infinite 
order from its holomorphic solutions.

\section{The functor $\IB:\ctTc\rightarrow\Ind(\ctBan)$}

For any object $E$ of $\ctTc$, we denote by $\ctb_{E}$ the set of 
absolutely convex bounded subsets of $E$ and by $\overline{\ctb}_{E}$ 
the set of closed absolutely convex bounded subsets of $E$.  If 
$B\in\ctb_{E}$, we denote $E_{B}$ the semi-normed space obtained by 
endowing the linear hull of $B$ in $E$ with the gauge semi-norm 
$p_{B}$ associated to $B$.

\begin{definition}
To define the functor
$$
\IB:\ctTc\to\Ind(\ctBan)
$$
we proceed as follows. For any object $E$ of $\ctTc$, we set
$$
\IB(E)=\findlim_{B\in\ctb_{E}}\widehat{E}_{B}
$$
where $\widehat{E}_{B}$ denotes as usual the completion of $E_{B}$.  
Consider a morphism $f:E\to F$ of $\ctTc$.  For any $B\in\ctb_{E}$, 
$f(B)\in\ctb_{F}$.  Hence, $f$ induces a morphism 
$\widehat{E}_{B}\to\widehat{F}_{f(B)}$.  This morphism being 
functorial in $B$, we obtain a morphism
$$
\findlim_{B\in\ctb_{E}}\widehat{E}_{B}\to
\findlim_{B\in\ctb_{E}}\widehat{F}_{f(B)}
$$
in $\Ind(\ctBan)$. We define
$$
\IB(f):\IB(E)\to\IB(F)
$$
by composing the preceding morphism with the canonical morphism
$$
\findlim_{B\in\ctb_{E}}\widehat{F}_{f(B)}\to
\findlim_{B\in\ctb_{F}}\widehat{F}_{B}.
$$
\end{definition}

\begin{remark}\label{rem:ibofban}
If $E$ is a Banach space, then
$$
\IB(E)\simeq\cindobj{E}.
$$
As a matter of fact, since any bounded subset of $E$ is included in a 
ball $b(\rho)$ centered at the origin, we have
$$
\IB(E)\simeq\findlim_{\rho>0}E_{b(\rho)}
$$
and the conclusion follows from the isomorphism $E_{b(\rho)}\simeq E$.
\end{remark}

\begin{lemma}\label{lem:isoboundedappli}
Let $E$ and $F$ be two objects of $\ctTc$. Then,
$$
\prolim_{B\in\ctb_{E}}\indlim_{B'\in\ctb_{F}}
\Hom[\ctTc](E_{B},F_{B'})\simeq B(E,F)
$$
where
$$
B(E,F)=\{f:E\to F:f\text{ linear},f(B)\text{ bounded in 
}F\text{ if } B\text{ bounded in }E\}.
$$
\end{lemma}

\begin{remark}\label{rem:bornhombound}
If $E$ and $F$ are objects of $\ctTc$, we have
$$
\Hom[\ctTc](E,F)\subset B(E,F).
$$
In general, this inclusion is strict but, as is well-known, it turns 
into an equality if $E$ is bornological (i.e.\ if any absolutely 
convex subset of $E$ that absorbs any bounded subset is a neighborhood 
of zero).
\end{remark}

\begin{proposition}\label{prp:borncomplhom}
Let $E$ and $F$ be two objects of $\ctTc$. If $E$ is bornological and 
$F$ complete, then
$$
\Hom[\Ind(\ctBan)](\IB(E),\IB(F))\simeq
\Hom[\ctTc](E,F).
$$
\end{proposition}

\begin{proof}
Since the inclusion $\overline{\ctb}_{F}\subset\ctb_{F}$ is cofinal, 
we have
\begin{align*}
\Hom[\Ind(\ctBan)](\IB(E),\IB(F))
  & \simeq\Hom[\Ind(\ctBan)]
    (\findlim_{B\in\ctb_{E}}\widehat{E}_{B},
    \findlim_{B'\in\overline{\ctb}_{F}}\widehat{F}_{B'})
\\
  & \simeq \prolim_{B\in\ctb_{E}}\indlim_{B'\in\overline{\ctb}_{F}}
    \Hom[\ctBan](\widehat{E}_{B},\widehat{F}_{B'}).
\end{align*}
Since $F$ is complete, $F_{B'}$ is a Banach space and
$$
\Hom[\ctBan](\widehat{E}_{B},\widehat{F}_{B'})\simeq
\Hom[\ctTc](E_{B},F_{B'}).
$$
It follows that
\begin{align*}
\Hom[\Ind(\ctBan)](\IB(E),\IB(F))
  & \simeq \prolim_{B\in\ctb_{E}}\indlim_{B'\in\ctb_{F}}
    \Hom[\ctTc](E_{B},F_{B'})
\\
  & \simeq B(E,F)
    \simeq \Hom[\ctTc](E,F)
\end{align*}
where the second isomorphism follows from 
Lemma~\ref{lem:isoboundedappli} and the last isomorphism from 
Remark~\ref{rem:bornhombound}.
\end{proof}

\begin{proposition}\label{prp:ilibadj}
Denote 
$$
\iolim:\Ind(\ctBan)\to\ctTc
$$
the functor defined by
$$
\iolim(\findlim_{i\in\cti}E_{i})=
\indlim_{i\in\cti}E_{i}.
$$
Let $E$ be an object of $\Ind(\ctBan)$ and let $F$ be a complete 
object of $\ctTc$.  Then,
$$
\Hom[\Ind(\ctBan)](E,\IB(F))\simeq\Hom[\ctTc](\iolim(E),F).
$$
\end{proposition}

\begin{proof}
Assuming $E\simeq\findlim_{i\in\cti}E_{i}$, we have
\begin{align*}
\Hom[\Ind(\ctBan)](E,\IB(F))
  & \simeq\prolim_{i\in\cti}\Hom[\Ind(\ctBan)](\cindobj{E_{i}},\IB(F))
\\
  & \simeq\prolim_{i\in\cti}\Hom[\Ind(\ctBan)](\IB(E_{i}),\IB(F))
\\
  & \simeq\prolim_{i\in\cti}\Hom[\ctTc](E_{i},F)
    \simeq\Hom[\ctTc](\iolim(E),F)
\end{align*}
where the second isomorphism follows from Remark~\ref{rem:ibofban} and 
the third from Proposition~\ref{prp:borncomplhom}.
\end{proof}

\begin{corollary}\label{cor:ibprolim}
Let $\cti$ be a small category.  For any functor 
$$
X:\cti^{\op}\to\ctTc
$$
such that $X(i)$ is complete for any $i\in\cti$, we have
$$
\IB(\prolim_{i\in\cti}X(i))\simeq\prolim_{i\in\cti}\IB(X(i)).
$$
\end{corollary}

\begin{proof}
For any object $E$ of $\Ind(\ctBan)$, we have
\begin{align*}
\Hom[\Ind(\ctBan)](E,\IB(\prolim_{i\in\cti}X(i)))
  &\simeq\Hom[\ctTc](\iolim(E),\prolim_{i\in\cti}X(i))
\\
  &\simeq\prolim_{i\in\cti}\Hom[\ctTc](\iolim(E),X(i))
\\
  &\simeq\prolim_{i\in\cti}\Hom[\Ind(\ctBan)](E,\IB(X(i)))
\end{align*}
where the first and last isomorphisms follow from 
Proposition~\ref{prp:ilibadj}.  The conclusion follows from the theory 
of representable functors.
\end{proof}

\begin{proposition}\label{prp:ibinjectacyclfr}
Assume that $(F_{n},f_{m,n})_{n\in\N}$ is an inductive system of 
Fr\'{e}chet spaces with injective transition morphisms and that
$$
\indlim_{n\in\N}F_{n}
$$
is complete.  Then, the canonical morphism
$$
\indlim_{n\in\N}\IB(F_{n})\to
\IB(\indlim_{n\in\N}F_{n})
$$
is an isomorphism.
\end{proposition}

\begin{proof}
Applying $\IB$ to the canonical morphisms
$$
r_{n}:F_{n}\to\indlim_{n\in\N}F_{n}
$$
and using the characterization of inductive limits, we get the 
canonical morphism
\begin{equation}\tag{*}
\indlim_{n\in\N}\IB(F_{n})\to
\IB(\indlim_{n\in\N}F_{n}).
\end{equation}
Let $B$ be a closed absolutely convex bounded subset of 
$\indlim_{n\in\N}F_{n}$.  It follows from e.g.~\cite[Chap.~IV, \S~19, 
5.(5) (p.~225)]{Kothe69} that, for some $n\in\N$, $B$ is the image of 
a closed absolutely convex bounded subset $B_{n}$ of $F_{n}$ by the 
canonical morphism $r_{n}$.  Since $r_{n}$ is injective, it induces 
the isomorphism of semi-normed spaces
$$
(F_{n})_{B_{n}}\isoto
(\indlim_{n\in\N}F_{n})_{B}.
$$
Hence, we get the isomorphism of Banach spaces
$$
(\widehat{\indlim_{n\in\N}F_{n}})_{B}\isoto
\widehat{(F_{n})}_{B_{n}}.
$$
Composing with the morphism
$$
\cindobj{\widehat{(F_{n})}_{B_{n}}}\to
\IB(F_{n})\to\indlim_{n\in\N}\IB(F_{n}),
$$
we get a canonical morphism
$$
\cindobj{(\widehat{\indlim_{n\in\N}F_{n}})_{B}}\to
\indlim_{n\in\N}\IB(F_{n}).
$$
Finally, using the characterization of inductive limits, we obtain a 
canonical morphism
$$
\IB(\indlim_{n\in\N}F_{n})=
\indlim_{B\in\ctb_{\indlim F_{n}}}
\cindobj{(\widehat{\indlim_{n\in\N}F_{n}})_{B}}
\to\indlim_{n\in\N}\IB(F_{n}).
$$
A direct computation shows that this morphism is a left and right 
inverse of (*).
\end{proof}

\begin{remark}
Note that, thanks to~\cite[Proposition~7.2]{Palamodov71} 
and~\cite[Corollary~7.2]{Palamodov71}, a countable filtering inductive 
system of Fr\'{e}chet spaces which is $\indlim$-acyclic in $\ctTc$ is 
essentially equivalent to an inductive system which satisfies the 
assumptions of the preceding proposition.  Hence, $\IB$ also commutes 
with the inductive limit functor in such a situation.
\end{remark}

\begin{definition}
Let $E$ and $F$ be two objects of $\ctTc$. As usual, we denote by 
$\Homtop[b](E,F)$ the vector space $\Hom[\ctTc](E,F)$ endowed with 
the system of semi-norms
$$
\{p_{B}:p\text{ continuous semi-norm of }F,\;B\text{ 
bounded subset of }E\}
$$
where
$$
p_{B}(f)=\sup_{e\in B}p(f(e)).
$$
\end{definition}

\begin{lemma}\label{lem:bornhomtopisoprolim}
Let $E$ and $F$ be two objects of $\ctTc$.  Assume $E$ is 
bornological.  Then,
$$
\Homtop[b](E,F)\simeq\prolim_{B\in\ctb_{E}}\Homtop[b](E_{B},F)
$$
in $\ctTc$. Assume moreover that $F$ is complete. Then,
$$
\Homtop[b](E,F)\simeq
\prolim_{B\in\ctb_{E}}\Homtop[b](\widehat{E}_{B},F)
$$
in $\ctTc$.
\end{lemma}

\begin{proof}
Keeping in mind the properties of bornological spaces, it is clear 
from the definition of $\Homtop[b](E,F)$ that
$$
\Homtop[b](E,F)\simeq
\prolim_{B\in\ctb_{E}}\Homtop[b](E_{B},F).
$$
Since any ball of $\widehat{E}_{B}$ is included in the closure of a 
semi-ball of $E_{B}$, any bounded subset of $\widehat{E}_{B}$ is 
included in the closure of a bounded subset of $E_{B}$. This property 
and the completeness of $F$ shows that
$$
\Homtop[b](E_{B},F)\simeq\Homtop[b](\widehat{E}_{B},F).
$$
Hence the conclusion.
\end{proof}

\begin{lemma}\label{lem:bancomplibhomtop}
If $E$ is a Banach space and if $F$ is a complete object of $\ctTc$, 
then
$$
\IB(\Homtop[b](E,F))\simeq
\Homtop(\cindobj{E},\IB(F)).
$$
\end{lemma}

\begin{proof}
For any $B'\in\ctb_{F}$, set
$$
B'_{b}=\{f\in\Hom[\ctTc](E,F):\norm{e}\leq 1\implies f(e)\in B'\}.
$$
Clearly, $B'_{b}$ belongs to $\ctb_{\Homtop[b](E,F)}$.
Moreover, if $B'$ is closed in $F$, then $B'_{b}$ is closed in 
$\Homtop[b](E,F)$ and one checks easily that
$$
(\Homtop[b](E,F))_{B'_{b}}\simeq\Homtop(E,F_{B'})
$$
as Banach spaces.  Hence, one has successively
\begin{align*}
\IB(\Homtop[b](E,F))
  & =\findlim_{B\in\overline{\ctb}_{\Homtop[b](E,F)}}
    (\Homtop[b](E,F))_{B}
   \simeq\findlim_{B'\in\overline{\ctb}_{F}}
    (\Homtop[b](E,F))_{B'_{b}}
\\
  & \simeq\findlim_{B'\in\overline{\ctb}_{F}}\Homtop(E,F_{B'})
    \simeq\Homtop(\cindobj{E},\IB(F))
\end{align*}
where the second isomorphism follows from the fact that the inclusion
$$
\{B'_{b}:B'\in\overline{\ctb}_{F}\}\subset
\overline{\ctb}_{\Homtop[b](E,F)}
$$
is cofinal.
\end{proof}

\begin{proposition}\label{prp:isoibhomtop}
Let $E$ and $F$ be two objects of $\ctTc$. Assume $E$ bornological 
and $F$ complete. Then,
$$
\IB(\Homtop[b](E,F))\simeq\Homtop(\IB(E),\IB(F)).
$$
\end{proposition}

\begin{proof}
We have successively
\begin{align}
\IB(\Homtop[b](E,F))
  & \simeq\IB(\prolim_{B\in\ctb_{E}}\Homtop[b](\widehat{E}_{B},F))
\tag{1}\\
  & \simeq\prolim_{B\in\ctb_{E}}\IB(\Homtop[b](\widehat{E}_{B},F))
\tag{2}\\
  & \simeq\prolim_{B\in\ctb_{E}}\Homtop(\cindobj{\widehat{E}_{B}},
    \IB(F))
\tag{3}\\
  & \simeq\Homtop(\IB(E),\IB(F)),
\notag
\end{align}
where the isomorphism (1) follows from 
Lemma~\ref{lem:bornhomtopisoprolim}, (2) from 
Corollary~\ref{cor:ibprolim} and (3) from 
Lemma~\ref{lem:bancomplibhomtop}.
\end{proof}

\begin{remark}
(1) Let $E$, $F$, $G$ be three objects of $\ctTc$. Recall that a 
bilinear application 
$$
b:E\times F\to G
$$
is continuous if and only if for any continuous semi-norm $r$ of $G$, 
there are continuous semi-norms $p$ and $q$ of $E$ and $F$ 
respectively such that
$$
r(b(x,y))\leq p(x)q(y).
$$

(2) Let $E$, $F$ be two objects of $\ctTc$ with $P$ and $Q$ as systems 
of semi-norms.  As usual, if $p\in P$ and $q\in Q$, we denote $p\tens 
q$ the semi-norm on $E\tens F$ defined by
$$
(p\tens q)(u)=
\inf_{u=\sum x_{i}\tens y_{i}}\sum_{i}p(x_{i})q(y_{i}).
$$
Recall that $E\tens_{\pi}F$ is the object of $\ctTc$ obtained by 
endowing $E\tens F$ with the system of semi-norms induced by 
$$
\{p\tens q:p\in P,q\in Q\}.
$$
From this definition, it follows immediately that any continuous 
bilinear map 
$$
b:E\times F\to G
$$
factors uniquely through a continuous linear map
$$
E\tens_{\pi}F\to G.
$$
Finally, recall that $E\tenstop_{\pi} F$ denotes the completion of 
$E\tens_{\pi}F$ and that $p\tenstop_{\pi} q$ is the semi-norm of 
$E\tenstop_{\pi} F$ induced by $p\tens q$.
\end{remark}

\begin{proposition}\label{prp:bilintenstpcom}
There is a canonical morphism
$$
\IB(E)\tenstop\IB(F)\to\IB(E\tens_{\pi}F).
$$
\end{proposition}

\begin{proof}
For $B\in\ctb_{E}$ and $B'\in\ctb_{F}$, denote $B\tens B'$ the 
absolutely convex hull of 
$$
\{b\tens b':b\in B,b'\in B'\}.
$$
This is clearly a bounded absolutely convex subset of $E\tens F$.  As 
a matter of fact,
$$
(p\tens q)(b\tens b')\leq p(b)q(b')\leq 
\sup_{b\in B}p(b)\sup_{b'\in B'}q(b').
$$
Moreover, we have a canonical linear map
$$
E_{B}\tens F_{B'}\to(E\tens_{\pi}F)_{B\tens B'}.
$$
This map is clearly continuous since $e\tens f\in B\tens B'$ when 
$e\in B$, $e'\in B'$. Applying the completion functor, we get a 
morphism
$$
\widehat{E}_{B}\tenstop\widehat{F}_{B'}\to
\widehat{(E\tens_{\pi}F)}_{B\tens B'}
$$
and hence a morphism
$$
\cproobj{\widehat{E}_{B}}\tenstop\cproobj{\widehat{F}_{B'}}
\simeq\IB(\widehat{E}_{B}\tenstop\widehat{F}_{B'})
\to\IB(E\tens_{\pi}F).
$$
Using the definition of inductive limits, we get a morphism
$$
\IB(E)\tenstop\IB(F)\simeq
\indlim_{B\in\ctb_{E}}\indlim_{B'\in\ctb_{F}}
\cproobj{\widehat{E}_{B}}\tenstop\cproobj{\widehat{F}_{B'}}
\to\IB(E\tens_{\pi}F).
$$
\end{proof}

\section{Some acyclicity results for $\Homtop$ and $\tenstop$ in 
$\Ind(\ctBan)$}

Hereafter, we denote as usual $c^{0}$ (resp.\ $l^{1}$) the Banach 
spaces formed by the sequences $x=(x_{n})_{n\in\N}$ of complex numbers 
that converge to $0$ (resp.\ that are summable); the norm being 
defined by
$$
\norm{x}_{c^{0}}=\sup_{n\in\N}|x_{n}|
\qquad\text{(resp.\ }
\norm{x}_{l^{1}}=\sum_{n=0}^{\infty}|x_{n}|\text{ )}.
$$
For any Banach space $X$, we also set for short 
$\Dual(X)=\Homtop(X,\C)$.

\begin{lemma}\label{lem:nucldecomp}
Let $X$, $Y$ be two Banach spaces and let $f:X\to Y$ be a nuclear map.  
Then, there is a continuous linear map $p:X\to c^{0}$ and a nuclear 
map $c:c^{0}\to Y$ making the diagram
$$
\xymatrix{
  & c^{0} \ar[rd]^-{c}
  &
\\
X \ar[ru]^-{p} \ar[rr]_-{f}
  &
  & Y
}
$$
commutative.
\end{lemma}

\begin{proof} 
Since $f:X\to Y$ is nuclear, there is a bounded sequence $x^{*}_{n}$ 
of $\Dual(X)$, a bounded sequence $y_{n}$ of $Y$ and a summable 
sequence $\lambda_{n}$ of complex numbers such that
$$
f(x)=\sum_{n=0}^{+\infty}\lambda_{n}\scal{x^{*}_{n},x}y_{n}
\qquad \forall x\in X.
$$
Since $\lambda_{n}$ is summable, one can find a sequence $r_{n}$ of 
non-zero complex numbers converging to zero such that 
$\lambda_{n}/r_{n}$ is still summable.  One checks easily that the 
maps $p:X\to c^{0}$ and $c:c^{0}\to Y$ defined by
$$
p(x)_{n}=r_{n}\scal{x^{*}_{n},x}
\qquad\text{and}\qquad
c(s)=\sum_{n=0}^{+\infty}
\frac{\lambda_{n}}{r_{n}}y_{n}s_{n}
$$
have the requested properties.
\end{proof}

\begin{definition}
A projective system $E:I^{\op}\to\ctBan$ where $I$ is a filtering 
ordered set is \emph{nuclear} if for any $i\in I$, there is $j\in I$, 
$j\geq i$ such that the transition morphism
$$
e_{i,j}:E_{j}\to E_{i}
$$
is nuclear.
\end{definition}

\begin{lemma}\label{lem:nuclprojsysczero}
Let $I$ be an infinite filtering ordered set and let 
$E:I^{\op}\to\ctBan$ be a nuclear projective system. Then, in 
$\Pro(\ctBan)$, we have
$$
\fprolim_{i\in I}E_{i}\simeq
\fprolim_{k\in K}X_{k}.
$$
where $X:K^{\op}\to\ctBan$ is a projective system with nuclear 
transition morphisms such that $X^{k}=c^{0}$ for any 
$k\in K$ and $\# K=\#I$.
\end{lemma}

\begin{proof}
Consider the set
$$
K=\{(i,j)\in I\times I:j\geq i,\; e_{i,j}:E_{j}\to E_{i}
\text{ nuclear}\}.
$$
The relation ``$\geq$'' defined by setting $(i',j')\geq (i,j)$ if 
$(i',j')=(i,j)$ or $i'\geq j$ turns $K$ into a filtering ordered set.  
By Lemma~\ref{lem:nucldecomp}, for any $k=(i,j)\in K$, we may choose a 
continuous linear map $p_{k}:E_{j}\to c^{0}$ and a nuclear map 
$c_{k}:c^{0}\to E_{i}$ making the diagram
$$
\xymatrix{
  & c^{0} \ar[rd]^-{c_{k}}
  &
\\
E_{j} \ar[ru]^-{p_{k}} \ar[rr]_-{e_{i,j}}
  &
  & E_{i}
}
$$
commutative. For any $k\in K$, we set $X_{k}=c^{0}$ and 
$x_{k,k}=\id_{X_{k}}$. If $k'=(i',j')>k=(i,j)$, we set
$$
x_{k,k'}=p_{k}\comp e_{j,i'}\comp c_{k'}:
X_{k'}\to X_{k}.
$$
The map $c_{k'}$ being nuclear, $x_{k,k'}$ is also nuclear. An easy 
computation shows that if $k<k'<k''$, then $x_{k,k'}\comp 
x_{k',k''}=x_{k,k''}$. Consider the functors
$$
\Phi:K\to I \qquad\text{and}\qquad
\Psi:K\to I
$$
defined by $\Phi((i,j))=i$ and $\Psi((i,j))=j$. They are clearly 
cofinal and if $k'\geq k$ in $K$, the diagrams
$$
\xymatrix{
X_{k'} \ar[r]^-{c_{k'}} \ar[d]_-{x_{k,k'}}
  & E_{\Phi(k')} \ar[d]^-{e_{\Phi(k),\Phi(k')}}
  & \qquad\qquad
  & E_{\Psi(k')} \ar[r]^-{p_{k'}}\ar[d]_-{e_{\Psi(k),\Psi(k')}}
  & X_{k'}\ar[d]^-{x_{k,k'}}
\\
X_{k} \ar[r]_-{c_{k}}
  & E_{\Phi(k)}
  & \qquad\qquad
  & E_{\Psi(k)} \ar[r]_-{p_{k}}
  & X_{k}
}
$$
are commutative. Hence, we get the two morphisms
$$
\fprolim_{k\in K}X_{k}\to\fprolim_{k\in K}E_{\Phi(k)}
\simeq\fprolim_{i\in I}E_{i}
\quad\text{and}\quad
\fprolim_{j\in I}E_{j}\simeq\fprolim_{k\in K}E_{\Psi(k)}
\to\fprolim_{k\in K}X_{k}.
$$
Since these morphisms are easily checked to be inverse one of each 
other, the proof is complete.
\end{proof}

\begin{remark}
Hereafter, as usual, we denote $e_{n}$ the element of $c^{0}$ defined 
by
$$
(e_{n})_{m}=\delta_{n,m}
$$
and we denote $e_{n}^{*}$ the element of $\Dual(c^{0})$ defined by
$$
\scal{e_{n}^{*},x}=x_{n}.
$$
\end{remark}

\begin{lemma}\label{lem:nuclczerotoban}
For any Banach space $Y$ and any nuclear map
$$
u:c^{0}\to Y
$$
the sequence $\norm{u(e_{n})}_{Y}$ is summable and for any $x\in 
c^{0}$, we have
$$
u(x)=\sum_{n=0}^{+\infty}\scal{e^{*}_{n},x}u(e_{n}).
$$
\end{lemma}

\begin{proof}
Since $u$ is nuclear, we can find a bounded sequence $x^{*}_{n}$ of 
$\Dual(c^{0})$, a bounded sequence $y_{n}$ of $Y$ and a summable 
sequence $\lambda_{n}$ of complex numbers such that
$$
u(x)=\sum_{n=0}^{+\infty}\lambda_{n}\scal{x^{*}_{n},x}y_{n}
$$
for any $x\in c^{0}$. Using the isomorphism $\Dual(c^{0})\simeq 
l^{1}$, we see that
\begin{equation}\tag{*}
\sum_{m=0}^{+\infty}|\scal{x_{n}^{*},e_{m}}|
=\norm{x_{n}^{*}}_{\Dual(c^{0})}.
\end{equation}
Therefore,
\begin{align*}
\sum_{m=0}^{M}\norm{u(e_{m})}
  &\leq \sum_{m=0}^{M}\sum_{n=0}^{+\infty}|\lambda_{n}|
    |\scal{x_{n}^{*},e_{m}}|\norm{y_{n}}_{Y}
\\
  &\leq \sum_{n=0}^{+\infty}|\lambda_{n}|
    \norm{x_{n}^{*}}_{\Dual(c^{0})}\norm{y_{n}}_{Y}
\\
  &\leq \left(\sum_{n=0}^{+\infty}|\lambda_{n}|\right)
    \sup_{n\in\N}\norm{x_{n}^{*}}_{\Dual(c^{0})}
    \sup_{n\in\N}\norm{y_{n}}_{Y}
\end{align*}
and the sequence $\norm{u(e_{n})}_{Y}$ is summable. Moreover,
\begin{align*}
u(x)
  &=\sum_{n=0}^{+\infty}\sum_{m=0}^{+\infty}\lambda_{n}
    \scal{x_{n}^{*},e_{m}}x_{m}y_{n}
\\
  &=\sum_{m=0}^{+\infty}x_{m}\left(\sum_{n=0}
    \lambda_{n}\scal{x_{n}^{*},e_{m}}y_{n}\right)
\\
  &=\sum_{m=0}^{+\infty}\scal{e_{m}^{*},x}u(e_{m})
\end{align*}
where the permutation of the sums is justified using (*).
\end{proof}

\begin{lemma}\label{lem:isotenstophomtopb}
Let $I$ be an infinite filtering ordered set and let 
$X:I^{\op}\to\ctBan$ be a nuclear projective system.  Assume $Y$ is a 
Fr\'{e}chet space.  Then, the morphisms
$$
\varphi_{i}:X_{i}\tenstop_{\pi} Y\to \Homtop[b](\Dual(X_{i}),Y)
$$
defined by setting
$$
\varphi_{i}(x\tenstop_{\pi} y)(x^{*})=\scal{x^{*},x}y
\qquad\forall x\in X_{i},\; y\in Y,\; x^{*}\in\Dual(X_{i})
$$
induce an isomorphism
$$
\fprolim_{i\in I}X_{i}\tenstop_{\pi} Y\simeq
\fprolim_{i\in I}\Homtop[b](\Dual(X_{i}),Y).
$$
In particular, for $Y=\C$, we have
$$
\fprolim_{i\in I}X_{i}\simeq
\fprolim_{i\in I}\Dual(\Dual(X_{i})).
$$
\end{lemma}

\begin{proof}
By Lemma~\ref{lem:nuclprojsysczero}, we may assume that $X_{i}=c^{0}$ 
for any $i\in I$ and that the transition morphisms
$$
x_{i,j}:X_{j}\to X_{i} \qquad(j>i)
$$
are nuclear.

One checks easily that $\varphi_{i}$ is a well-defined continuous map.  
By Lemma~\ref{lem:nuclczerotoban}, we know that the sequence 
$(\norm{x_{i,j}(e_{n})}_{X_{i}})_{n\in\N}$ is summable and that
$$
x_{i,j}(c)=\sum_{n=0}^{+\infty}\scal{e^{*}_{n},c}x_{i,j}(e_{n})
\qquad\forall c\in X_{j}.
$$
Therefore, we may define a continuous linear map
$$
\psi_{i,j}:\Homtop[b](\Dual(X_{j}),Y)\to
X_{i}\tenstop_{\pi}Y
$$
by setting
$$
\psi_{i,j}(h)=\sum_{n=0}^{+\infty}x_{i,j}(e_{n})
\tenstop_{\pi}h(e_{n}^{*}).
$$
One sees easily that the morphisms $\varphi_{i}$ and $\psi_{i,j}$ 
induce morphisms of pro-objects
$$
\fprolim_{i\in I}X_{i}\tenstop_{\pi}Y\to
\fprolim_{i\in I}\Homtop[b](\Dual(X_{i}),Y)
$$
and
$$
\fprolim_{i\in I}\Homtop[b](\Dual(X_{i}),Y)\to
\fprolim_{i\in I}X_{i}\tenstop_{\pi}Y.
$$
A direct computation shows that these morphisms are inverse one of 
each other.
\end{proof}

\begin{definition}
We say that a filtering projective system $E:I^{\op}\to\ctTc$ 
satisfies \emph{Condition~ML} if for any $i\in I$, any semi-norm $p$ 
of $E_{i}$ and any $\epsilon>0$, there is $i'\geq i$ such that
$$
e_{i,i'}(E_{i})\subset b_{p}(\epsilon)+e_{i,i''}(E_{i''})
\qquad\forall i''\geq i'
$$
where $b_{p}(\epsilon)$ denotes as usual the semi-ball of radius 
$\epsilon$ and center $0$ associated to the semi-norm $p$.
\end{definition}

\begin{remark}\label{rem:limacyclic}
By~\cite[Proposition~1.2.9]{Prosmans99a} (which is a direct 
consequence of~\cite[Theorem~5.6]{Prosmans99}), a countable filtering 
projective system of Fr\'{e}chet spaces is $\prolim$-acyclic in 
$\ctTc$ if and only if it satisfies Condition~ML.
\end{remark}

\begin{lemma}\label{lem:tenstopmitlef}
Let $E:I^{\op}\to\ctTc$ and $F:J^{\op}\to\ctTc$ be two filtering 
projective systems.  If $E$ and $F$ satisfy Condition~ML, then the 
projective system
$$
E\tenstop_{\pi}F:(I\times J)^{\op}\to\ctTc
$$
defined by
$$
(E\tenstop_{\pi}F)(i,j)=E_{i}\tenstop_{\pi}F_{j}
$$
satisfies Condition~ML.
\end{lemma}

\begin{proof}
Let $(i,j)\in I\times J$ and let $p\tenstop_{\pi}q$ be a semi-norm of 
$E_{i}\tenstop_{\pi}F_{j}$.  It follows from our assumptions, that 
there is $i'\geq i$ and $j'\geq j$ such that
\begin{equation}\tag{*}
e_{i,i'}(E_{i'})\subset b_{p}(1)+e_{i,i''}(E_{i''})
\qquad \forall i''\geq i'
\end{equation}
and
\begin{equation}\tag{**}
f_{j,j'}(F_{j'})\subset b_{q}(1)+f_{j,j''}(F_{j''})
\qquad\forall j''\geq j'.
\end{equation}
Fix $(i'',j'')\geq (i',j')$.  Since the maps $e_{i,i'}$, $f_{j,j'}$ 
and $e_{i,i''}$ are continuous, we can find a semi-norm $p'$ of 
$E_{i'}$, a semi-norm $q'$ of $F_{j'}$ and a semi-norm $p''$ of 
$E_{i''}$ such that
$$
p\comp e_{i,i'}\leq p',\qquad q\comp f_{j,j'}\leq q'
\qquad\text{and}\qquad
p\comp e_{i,i''}\leq p''.
$$
Consider $\epsilon>0$ and let $z'$ be an element of $E_{i'}\tens_{\pi} 
F_{j'}$ of the type $x'\tens_{\pi}y'$ where $x'\in E_{i'}$, $y'\in 
F_{j'}$.  Using (*) and (**) above, we obtain $x''\in E_{i''}$ 
and $y''\in F_{j''}$ such that
$$
p(e_{i,i'}(x')-e_{i,i''}(x''))\leq \frac{\epsilon}{2(1+q'(y'))}
$$
and
$$
q(f_{j,j'}(y')-f_{j,j''}(y''))\leq \frac{\epsilon}{2(1+p''(x''))}.
$$
For $z''=x''\tens_{\pi}y''\in E_{i''}\tens_{\pi}F_{j''}$, we get
\begin{align*}
&(p\tens_{\pi}q)\left((e_{i,i'}\tens_{\pi}f_{j,j'})(z')-
  (e_{i,i''}\tens_{\pi}f_{j,j''})(z'')\right)
\\
&\quad =(p\tens_{\pi}q)
    \left((e_{i,i'}(x')-e_{i,i''}(x''))\tens_{\pi}f_{j,j'}(y')
    \right.
\\
&\qquad\qquad
    \left. 
      +e_{i,i''}(x'')\tens_{\pi}(f_{j,j'}(y')-f_{j,j''}(y''))
    \right)
\\
&\quad \leq p(e_{i,i'}(x')-e_{i,i''}(x''))q(f_{j,j'}(y'))+
    p(e_{i,i''}(x''))q(f_{j,j'}(y')-f_{j,j''}(y''))
    \leq \epsilon.
\end{align*}
Since any element of $E_{i}\tens_{\pi}F_{j}$ is a finite sum of 
elements of the type considered above, we see that for any 
$\epsilon>0$,
$$
(e_{i,i'}\tens_{\pi}f_{j,j'})(E_{i'}\tens_{\pi}F_{j'})\subset 
b_{p\tens_{\pi}q}(\epsilon)+
(e_{i,i''}\tens_{\pi}f_{j,j''})(E_{i''}\tens_{\pi}F_{j''}).
$$
The conclusion follows directly since $E_{i'}\tens_{\pi}F_{j'}$ is 
dense in $E_{i'}\tenstop_{\pi}F_{j'}$.
\end{proof}

\begin{remark}
Let $E$ be an object of $\ctTc$.  Recall that $E$ is of type FN if it 
is a nuclear Fr\'{e}chet space and that $E$ is of type DFN if it is 
isomorphic to the strong dual of a nuclear Fr\'{e}chet space.
\end{remark}

\begin{lemma}\label{lem:fnstruct}
Assume $X$ is a FN space.  Then, there is a projective system 
$$
(X_{n},x_{n,m})_{n\in\N}
$$
of Banach spaces such that
\begin{itemize}
\item[(a)]
there is an isomorphism
$$
X\simeq\prolim_{n\in\N}X_{n};
$$

\item[(b)]
for $m>n$, the transition map
$$
x_{n,m}:X_{m}\to X_{n}
$$
is nuclear and has a dense range;

\item[(c)]
there is an isomorphism
$$
\SDual(X)\simeq\indlim_{n\in\N}\Dual(X_{n})
$$
where $\SDual(X)$ denotes the strong dual of $X$;

\item[(d)]
for $m>n$, the transition map
$$
\Dual(x_{n,m}):\Dual(X_{n})\to\Dual(X_{m})
$$
is nuclear and injective.
\end{itemize}
\end{lemma}

\begin{proof}
Since $X$ is a FN space, there is a cofinal increasing sequence 
$(p_{n})_{n\in\N}$ of continuous semi-norms of $X$ such that the 
canonical map
$$
X_{p_{n+1}}\to X_{p_{n}}
$$
is nuclear. For such a sequence, the canonical map
$$
\widehat{X}_{p_{n+1}}\to \widehat{X}_{p_{n}}
$$
is also nuclear and has a dense range.  Moreover, it is well-known 
(see~e.g.~\cite[Chap.~IV, \S~19, 9.(1) (p.~231)]{Kothe69}) that 
$$
X\simeq\prolim_{n\in\N}\widehat{X}_{p_{n}}.
$$

Clearly,
$$
\IDual(X)\simeq\indlim_{n\in\N}\Dual(X_{p_{n}})
\simeq\indlim_{n\in\N}\Dual(\widehat{X}_{p_{n}})
$$
where $\IDual(X)$ is the inductive dual of $X$.

Recall that an absolutely convex subset $V$ is a neighborhood of $0$ 
in $\IDual(X)$ if it absorbs any equicontinuous subset of $X'$.  
Hence, it is clear that a neighborhood of $0$ in $\SDual(X)$ is a 
neighborhood of $0$ in $\IDual(X)$.  We know that $X$ is reflexive 
(see~e.g.~\cite[\S~5.3.2 (p.~93)]{Pietsch72}).  Hence, $\SDual(X)$ is 
bornological (see~e.g.~\cite[Chap.~VI, \S~29, 4.(4) 
(p.~400)]{Kothe69}).  The space $X$ being itself bornological, the 
bounded subsets of $\SDual(X)$ are equicontinuous.  So, any 
neighborhood of $0$ in $\IDual(X)$ is a neighborhood of $0$ in 
$\SDual(X)$ and $\IDual(X)\simeq\SDual(X)$.

Since (d) follows directly from (b), the proof is complete.
\end{proof}

\begin{proposition}\label{prp:dfnfribhomacycl}
Assume $E$ is a DFN space and $F$ is a Fr\'{e}chet space. Then, the 
canonical morphism
$$
\Hom(\IB(E),\IB(F))\to\RHom(\IB(E),\IB(F))
$$
is an isomorphism.
\end{proposition}

\begin{proof}
Since $E$ is a DFN space, there is a FN space $X$ such 
that
$$
E\simeq\SDual(X).
$$
Let $(X_{n},x_{n,m})$ be a projective system of the kind considered 
in Lemma~\ref{lem:fnstruct}. We have
$$
E\simeq\SDual(X)\simeq\indlim_{n\in\N}\Dual(X_{n}).
$$
Since the transition morphisms
$$
\Dual(x_{n,m}):\Dual(X_{n})\to\Dual(X_{m}) \qquad(m>n)
$$
are injective and $E$ is complete, 
Proposition~\ref{prp:ibinjectacyclfr} and Remark~\ref{rem:ibofban} 
show that
$$
\IB(E)\simeq\indlim_{n\in\N}\cindobj{\Dual(X_{n})}.
$$
Using Lemma~\ref{lem:nuclprojsysczero}, we find a nuclear projective 
system $(Y_{n},y_{n,m})$ with $Y_{n}=c^{0}$ such that
$$
\fprolim_{n\in\N}X_{n}\simeq\fprolim_{n\in\N}Y_{n}.
$$
It follows that 
$$
\IB(E)\simeq\indlim_{n\in\N}\cindobj{\Dual(Y_{n})}.
$$
Hence, we have successively
\begin{align}
\RHom(\IB(E),\IB(F))
  &\simeq\RHom(\lindlim_{n\in\N}\cindobj{\Dual(Y_{n})},\IB(F))
\tag{1}\\
  &\simeq\rprolim_{n\in\N}\RHom(\cindobj{\Dual(Y_{n})},\IB(F))
\tag{2}\\
  &\simeq\rprolim_{n\in\N}\Hom(\cindobj{\Dual(Y_{n})},\IB(F))
\tag{3}\\
  &\simeq\rprolim_{n\in\N}\Hom(\IB(\Dual(Y_{n})),\IB(F))
\tag{4}\\
&\simeq\rprolim_{n\in\N}\Hom[\ctTc](\Dual(Y_{n}),F) \tag{5}\end{align} 
where the isomorphism (1) follows from the fact that filtering 
inductive limits are exact in $\Ind(\ctBan)$, (2) follows 
from~\cite[Proposition~3.6.3]{Prosmans98a}, (3) follows from the fact 
that 
$
\cindobj{\Dual(Y_{n})}
\simeq
\cindobj{\Dual(c^{0})}\simeq\cindobj{l^{1}}
$ 
is projective in $\Ind(\ctBan)$, (4) follows from 
Remark~\ref{rem:ibofban} and (5) follows from 
Proposition~\ref{prp:borncomplhom}.  By 
Lemma~\ref{lem:isotenstophomtopb}, we have the isomorphism
$$
\fprolim_{n\in\N}(Y_{n}\tenstop_{\pi}F)\simeq
\fprolim_{n\in\N}\Homtop[b](\Dual(Y_{n}),F).
$$
Forgetting the topologies and applying the derived projective limit 
functor for pro-objects (see~\cite{Prosmans98a}), we obtain the 
isomorphism
$$
\rprolim_{n\in\N}(Y_{n}\tenstop_{\pi}F)\simeq
\rprolim_{n\in\N}\Hom[\ctTc](\Dual(Y_{n}),F).
$$
Since $(X_{n},x_{n,m})_{n\in\N}$ satisfies Condition~ML, it is 
$\prolim$-acyclic in $\ctTc$ (see Remark~\ref{rem:limacyclic}).  It 
follows that $(Y_{n},y_{n,m})_{n\in\N}$ is also $\prolim$-acyclic in 
$\ctTc$ and, hence, satisfies Condition~ML. Using 
Lemma~\ref{lem:tenstopmitlef}, we see that 
$\rprolim_{n\in\N}(Y_{n}\tenstop_{\pi}F)$ is concentrated in degree 
$0$.  It follows that the projective system
$$
\left(\Hom[\ctTc](\Dual(Y_{n}),F)\right)_{n\in\N}
$$
is $\prolim$-acyclic and the conclusion follows.
\end{proof}

\begin{theorem}\label{thm:dfnfrhomtopacyc}
Assume $E$ is a DFN space and $F$ is a Fr\'{e}chet space.  Then, the 
canonical morphism
$$
\Homtop(\IB(E),\IB(F))\to\RHomtop(\IB(E),\IB(F))
$$
is an isomorphism.
\end{theorem}

\begin{proof}
It is sufficient to show that
$$
LH^{k}(\RHomtop(\IB(E),\IB(F)))\simeq 0
$$
for $k>0$. This will be the case if 
$$
\Hom(\cindobj{l^{1}(I)},\RHomtop(\IB(E),\IB(F)))
$$
is concentrated in degree $0$ for any set $I$.

Let $I$ be an arbitrary set.  Since $\cindobj{l^{1}(I)}$ is a 
projective object of $\Ind(\ctBan)$ and since $F$ is complete, we have
\begin{align*}
\RHomtop(\cindobj{l^{1}(I)},\IB(F))
  & \simeq\Homtop(\IB(l^{1}(I)),\IB(F))
\\
  & \simeq\IB(\Homtop[b](l^{1}(I),F))
\\
  & \simeq\IB(l^{\infty}(I,F))
\end{align*}
where $l^{\infty}(I,F)$ is the Fr\'{e}chet space formed by the bounded 
families $(x_{i})_{i\in I}$ of $F$ (a fundamental system of semi-norms 
being given by 
$$
\{p_{I}:p \text{ continuous semi-norm of }F\}
$$
where $p_{I}((x_{i})_{i\in I})=\sup_{i\in I}p(x_{i})$).  Therefore, we 
have the chain of isomorphisms
\begin{align}
\Hom(\cindobj{l^{1}(I)},\RHomtop(\IB(E),\IB(F)))
  & \simeq\RHom(\cindobj{l^{1}(I)},\RHomtop(\IB(E),\IB(F)))
\notag\\
  & \simeq\RHom(\cindobj{l^{1}(I)}\ltenstop\IB(E),\IB(F))
\notag\\
  & \simeq\RHom(\IB(E)\ltenstop\cindobj{l^{1}(I)},\IB(F))
\notag\\
  & \simeq\RHom(\IB(E),\RHomtop(\cindobj{l^{1}(I)},\IB(F)))
\notag\\
  & \simeq\RHom(\IB(E),\IB(l^{\infty}(I,F)))
\notag
\end{align}
and the conclusion follows from Proposition~\ref{prp:dfnfribhomacycl}.
\end{proof}

\begin{lemma}\label{lem:nucldecompl1}
Let $X$, $Y$ be two Banach spaces and let $f:X\to Y$ be a nuclear map.  
Then, there is a nuclear map $p:X\to l^{1}$ and a continuous linear 
map $c:l^{1}\to Y$ making the diagram
$$
\xymatrix{
  & l^{1} \ar[rd]^-{c}
  &
\\
X \ar[ru]^-{p} \ar[rr]_-{f}
  &
  & Y
}
$$
commutative.
\end{lemma}

\begin{proof}
Work as for Lemma~\ref{lem:nucldecomp}.
\end{proof}

\begin{definition}
An inductive system $E:I\to\ctBan$ where $I$ is a filtering ordered 
set is \emph{nuclear} if for any $i\in I$, there is $j\in I$, $j\geq 
i$ such that the transition morphism
$$
e_{j,i}:E_{i}\to E_{j}
$$
is nuclear. An object of $\Ind(\ctBan)$ is \emph{nuclear} if it 
corresponds to a nuclear inductive system.
\end{definition}

\begin{remark}
Working as in the proof of Proposition~\ref{prp:dfnfribhomacycl}, we 
see easily that $\IB(E)$ is nuclear if $E$ is a DFN space.
\end{remark}

\begin{lemma}\label{lem:nuclindsysl1}
Let $I$ be an infinite filtering ordered set and let 
$E:I\to\ctBan$ be a nuclear inductive system. Then,
$$
\findlim_{i\in I}E_{i}\simeq
\findlim_{k\in K}X_{k}.
$$
where $X:K\to\ctBan$ is an inductive system with nuclear transition 
morphisms such that $X^{k}=l^{1}$ for any $k\in K$ and $\# K=\#I$.
\end{lemma}

\begin{proof}
Work as for Lemma~\ref{lem:nuclprojsysczero} using 
Lemma~\ref{lem:nucldecompl1}.
\end{proof}

\begin{lemma}\label{lem:tenstopacyclic}
Let $I$ be a filtering ordered set.  For any 
$F\in D^{-}(\Ind(\ctBan))$ and any $E\in D^{-}(\Ind(\ctBan)^{I})$,
we have
$$
(\indlim_{i\in I}E_{i})\ltenstop F\simeq
\indlim_{i\in I}(E_{i}\ltenstop F).
$$
\end{lemma}

\begin{proof}
If $P_{\cdot}$ is a projective resolution of $F$, we have successively
$$
(\indlim_{i\in I}E_{i})\ltenstop F\simeq
(\indlim_{i\in I}E_{i})\tenstop P_{\cdot}\simeq
\indlim_{i\in I}(E_{i}\tenstop P_{\cdot})\simeq
\indlim_{i\in I}(E_{i}\ltenstop F).
$$
\end{proof}

\begin{proposition}\label{prp:dfntensacyclic}
Let $E$ and $F$ be objects of $\Ind(\ctBan)$. Assume $E$ is nuclear. 
Then,
$$
E\ltenstop F\simeq E\tenstop F.
$$
\end{proposition}

\begin{proof}
Using Lemma~\ref{lem:nuclindsysl1}, we may assume that
$$
E=\findlim_{i\in I}X_{i}
$$
where $X:I\to\ctBan$ is a filtering inductive system with 
$X_{i}=l^{1}$, the transition morphisms
$$
x_{j,i}:X_{i}\to X_{j}
$$
being nuclear.  We may also assume that
$$
F=\findlim_{j\in J}Y_{j}
$$
where $Y:J\to\ctBan$ is a filtering inductive system. Then, we have
\begin{align}
E\ltenstop F
  &\simeq (\findlim_{i\in I}X_{i})\ltenstop(\findlim_{j\in J}Y_{j})
\notag\\
  &\simeq \indlim_{i\in I}\indlim_{j\in J}
      \cindobj{X_{i}}\ltenstop\cindobj{Y_{j}}
\tag{1}\\
  &\simeq \indlim_{i\in I}\indlim_{j\in J}
      \cindobj{X_{i}}\tenstop\cindobj{Y_{j}}
\tag{2}\\
  &\simeq (\findlim_{i\in I}X_{i})\tenstop(\findlim_{j\in J}Y_{j})
\notag\\
  &\simeq E\tenstop F
\notag
\end{align}
where the isomorphism (1) follows from Lemma~\ref{lem:tenstopacyclic} 
and (2) from the fact that 
$\cindobj{X_{i}}\simeq\cindobj{l^{1}}$ is projective in 
$\Ind(\ctBan)$.
\end{proof}

\section{A topological version of Cartan's Theorem B}

\begin{proposition}\label{prp:ibsheaf}
Let $X$ be a topological space with a countable basis.  If $F$ is a 
presheaf of Fr\'{e}chet spaces on $X$ which is a sheaf of vector 
spaces, then
$$
U\mapsto\IB(F(U))
\qquad(U\text{ open of }X)
$$
is a sheaf with values in $\Ind(\ctBan)$.
\end{proposition}

\begin{proof}
Let $U$ be an open subset of $X$ and let $\ctu$ be an open covering 
of $U$. Consider the sequence
\begin{equation}\tag{*}
0\to F(U)\to[\alpha]\prod_{V \in\ctu}F(V)
\to[\beta]\prod_{V,W\in\ctu}F(V\inter W)
\end{equation}
where $\alpha$ and $\beta$ are the continuous applications defined by
$$
p_{V}\comp\alpha=r_{V,U}
\qquad\text{and}\qquad
p_{V,W}\comp\beta=r_{V\inter W,V}\comp p_{V}-r_{V\inter W,W}\comp p_{W}
$$
where $p_{V}$ and $p_{V,W}$ are the canonical projections and 
$r_{V,U}$ is the restriction map. Since $F$ is a sheaf of vector 
spaces, this sequence is algebraically exact. Let us show that it is 
strictly exact.

(1) If $\ctu$ is countable, $F(U)$, $\prod_{V\in\ctu}F(V)$ and 
$\prod_{V,W\in\ctu}F(V\inter W)$ are Fr\'{e}chet spaces. Then, by the 
homomorphism theorem, the sequence (*) is strictly exact.

(2) Assume that $\ctu$ is not countable. Since $X$ has a countable 
basis, there is a countable set $A$ of open subsets of $X$ such that 
for any open $V$ of $X$,
$$
V=\Union_{k\in\N}U_{k},\qquad U_{k}\in A.
$$
Then, consider the countable set
$$
\ctv=\{V'\in A:\;\exists\,V\in\ctu\text{ such that } V'\subset V\}.
$$
For any $U'\in\ctu$, we may assume that $U'=\Union_{k\in\N}U'_{k}$, 
with $U'_{k}\in \ctv$.  It follows that $\ctv$ covers any $U'$ in 
$\ctu$ and therefore is a covering of $U$.  Hence, by (1), the sequence
$$
0\to F(U)\to[\alpha']\prod_{V'\in\ctv}F(V')
\to[\beta']\prod_{V',W'\in\ctv}F(V'\inter W')
$$
is strictly exact.  Now, consider a map $f:\ctv\to\ctu$ such that 
$V'\subset f(V')$ for any $V'\in\ctv$.  Then, consider the commutative 
diagram
$$
\xymatrix{
0 \ar[r]
  & F(U) \ar[r]^-{\alpha}\ar[d]^-{\id}
  & \prod_{V\in\ctu}F(V) \ar[r]^-{\beta}\ar[d]^-{\gamma}
  & \prod_{V,W\in\ctu}F(V\inter W) \ar[d]^-{\delta}
\\
0 \ar[r]
  & F(U) \ar[r]_-{\alpha'}
  & \prod_{V'\in\ctv}F(V') \ar[r]_-{\beta'}
  & \prod_{V',W'\in\ctv}F(V'\inter W') 
}
$$
where $\gamma$ and $\delta$ are respectively defined by
$$
p_{V'}\comp\gamma=r_{V',f(V')}\comp p_{f(V')}
$$
and
$$
p_{V',W'}\comp\delta=r_{V'\inter W',f(V')\inter f(W')}\comp 
p_{f(V'),f(W')}.
$$
To prove that the sequence (*) is strictly exact, it is sufficient to 
establish that $\alpha$ is a kernel of $\beta$. Let 
$h:X\to\prod_{V\in\ctu}F(V)$ be a morphism of $\ctTc$ such that 
$\beta\comp h=0$. Since $\beta'\comp\gamma\comp 
h=\delta\comp\beta\comp h=0$ and since $\alpha'$ is a kernel of 
$\beta'$, there is a unique morphism $h':X\to F(U)$ such that 
$\alpha'\comp h'=\gamma\comp h$. Set
$$
h''=h-\alpha\comp h'.
$$
We clearly have $\gamma\comp h''=0$ and $\beta\comp h''=0$. Fix 
$V\in\ctu$. For any $V'\in\ctv$ such that $V'\subset V$, we have
$$
0=p_{V,f(V')}\comp\beta\comp h''=
r_{V\inter f(V'),V}\comp p_{V}\comp h''-
r_{V\inter f(V'),f(V')}\comp p_{f(V')}\comp h''.
$$
It follows that
\begin{align*}
r_{V',V}\comp p_{V}\comp h''
  &=r_{V',V\inter f(V')}\comp
    r_{V\inter f(V'),V}\comp p_{V}\comp h''
\\
  &=r_{V',V\inter f(V')}\comp
    r_{V\inter f(V'),f(V')}\comp p_{f(V')}\comp h''
\\
  &=r_{V',f(V')}\comp p_{f(V')}\comp h''
\\
  &=p_{V'}\comp\gamma\comp h''
\\
  &=0.
\end{align*}
Since $\{V'\in\ctv:V'\subset V\}$ is a covering of $V$ and since $F$ 
is a sheaf of vector spaces, we get
$$
p_{V}\comp h''=0
\qquad\forall V\in\ctu.
$$
It follows that $h''=0$ and that $h=\alpha\comp h'$. Since $\alpha$ 
is injective, $h'$ is the unique morphism of $\ctTc$ such that 
$h=\alpha\comp h'$. Therefore, $\alpha$ is a kernel of $\beta$ and the 
sequence (*) is strictly exact.

Finally, since the functor $\IB$ preserves projective limits of 
complete objects of $\ctTc$ (see~Corollary~\ref{cor:ibprolim}), the 
sequence
$$  
0\to \IB(F(U))\to[\IB(\alpha)]\prod_{V \in\ctu}\IB(F(V))
\to[\IB(\beta)]\prod_{V,W\in\ctu}\IB(F(V\inter W))
$$
is strictly exact in $\Ind(\ctBan)$ .  Hence, the conclusion.
\end{proof}

\begin{definition}
For short, we denote $\IB(F)$ the sheaf with values in $\Ind(\ctBan)$ 
associated to a presheaf $F$ of the kind considered in 
Proposition~\ref{prp:ibsheaf}.
\end{definition}

Hereafter, $X$ will denote a complex analytic manifold of complex 
dimension $d_{X}$.  We denote $\O_{X}$ the sheaf of holomorphic 
functions on $X$.  Recall that for any open subset $U$ of $X$, 
$\O_{X}(U)$ has a canonical structure of FN space.  Recall moreover 
that if $V$ is a relatively compact open subset of $U$ the restriction 
morphism
$$
\O_{X}(U)\to\O_{X}(V)
$$
is nuclear. In particular, if $K$ is a compact subset of $X$, then
$$
\O_{X}(K)\simeq
\indlim_{\substack{U\supset K\\ U\text{open}}}\O_{X}(U)
$$
topologized as an inductive limit is a DFN space.

\begin{proposition}\label{prp:ibisosect}
For any compact subset $K$ of $X$, we have
$$
\Sect(K,\IB(\O_{X}))\simeq\IB(\O_{X}(K)).
$$
\end{proposition}

\begin{proof}
We know that $K$ has a fundamental system $(U_{n})_{n\in\N}$ of 
relatively compact open neighborhoods such that
$$
\overline{U}_{n+1}\subset U_{n}
$$
for any $n\in\N$.  Replacing, if necessary, $U_{n}$ by the union of 
those of its connected components which meet $K$, we may even assume 
that any connected component of $U_{n}$ meets $K$.  In this case, it 
follows from the principle of unique continuation that the restriction
$$
\O_{X}(U_{n})\to\O_{X}(U_{n+1})
$$
is injective. Moreover, by cofinality,
$$
\O_{X}(K)\simeq\indlim_{n\in\N}\O_{X}(U_{n}).
$$
Hence, by Proposition~\ref{prp:ibinjectacyclfr}, it follows 
that
$$
\IB(\indlim_{n\in\N}\O_{X}(U_{n}))\simeq
\indlim_{n\in\N}\IB(\O_{X}(U_{n})).
$$
Since $K$ is a taut subspace of $X$, a cofinality argument shows that
$$
\Sect(K,\IB(\O_{X}))\simeq
\indlim_{n\in\N}\Sect(U_{n},\IB(\O_{X}))
$$
and the conclusion follows.
\end{proof}

Hereafter, we denote $\ctc_{\infty,X}$ the sheaf of rings formed by 
functions of class $C_{\infty}$.  More generally, we denote 
$\ctc^{(p,q)}_{\infty,X}$ the sheaf of differential forms of class 
$C_{\infty}$ and of bitype $(p,q)$.  Recall that for any open subset 
$U$ of $X$, $\ctc^{(p,q)}_{\infty,X}(U)$ has a canonical structure of 
FN space.  Since the conditions of Proposition~\ref{prp:ibsheaf} are 
satisfied, $\IB(\ctc^{(p,q)}_{\infty,X})$ is a sheaf with values in 
$\Ind(\ctBan)$.

\begin{proposition}\label{prp:cinftyacyclic}
The sheaf $\IB(\ctc^{(p,q)}_{\infty,X})$ is 
$\Sect(U,\cdot)$-acyclic for any open subset $U$ of $X$.
\end{proposition}

\begin{proof}
For any object $E$ of $\Ind(\ctBan)$, denote 
$$
h_{E}:\Ind(\ctBan)\to\ctAb
$$
the functor defined by setting $h_{E}(F)=\Hom(E,F)$.  Using the 
techniques developed in~\cite{Schneiders99}, one shows easily that
$$
\Hom(P,\RSect(U,\IB(\ctc^{(p,q)}_{\infty,X})))
\simeq\RSect(U,h_{P}(\IB(\ctc^{(p,q)}_{\infty,X})))
$$
for any projective object $P$ of $\Ind(\ctBan)$.  Therefore, the 
result will be true if the sheaf of abelian groups 
$h_{P}(\IB(\ctc^{(p,q)}_{\infty,X}))$ is soft.  This follows 
from the fact that it has clearly a canonical structure of 
$\ctc_{\infty,X}$-module.
\end{proof}

\begin{theorem}\label{thm:ibsectacyclic}
If $U$ is an open subset of $X$ such that
$$
H^{k}(U,\O_{X})\simeq 0 \qquad(k>0)
$$
algebraically, then
$$
\RSect(U,\IB(\O_{X}))\simeq\IB(\O_{X}(U)).
$$
\end{theorem}

\begin{proof}
As is well-known, since $\ctc^{(p,q)}_{\infty,U}$ is a soft sheaf, the 
Dolbeault complex
$$
0\to\ctc^{(0,0)}_{\infty,X}
\to[\overline{\partial}]\ctc^{(0,1)}_{\infty,X}
\cdots\to[\overline{\partial}]\ctc^{(0,n)}_{\infty,X}\to 0
$$
is a $\Sect(U,\cdot)$-acyclic resolution of $\O_{X}$. Therefore,  
$\RSect(U,\O_{X})$ is given by the complex
$$
0\to\Sect(U,\ctc^{(0,0)}_{\infty,X})
\to[\overline{\partial}]\Sect(U,\ctc^{(0,1)}_{\infty,X})
\cdots\to[\overline{\partial}]\Sect(U,\ctc^{(0,n)}_{\infty,X})\to 0.
$$
Moreover, since $H^{k}(U,\O_{X})\simeq 0$ for $k>0$, the 
sequence
$$
0\to\Sect(U,\O_{X})\to\Sect(U,\ctc^{(0,0)}_{\infty,X})
\to[\overline{\partial}]\Sect(U,\ctc^{(0,1)}_{\infty,X})
\cdots\to[\overline{\partial}]\Sect(U,\ctc^{(0,n)}_{\infty,X})\to 0
$$
is algebraically exact.  Since $\O_{X}(U)$ and 
$\ctc^{(p,q)}_{\infty,X}(U)$ are FN spaces, the last sequence is 
strictly exact in $\ctTc$.  
Using~\cite[Proposition~3.2.26]{Schneiders99}, one sees easily that 
the sequence
\begin{equation}\tag{*}
0\to\Sect(U,\IB(\O_{X}))\to
\Sect(U,\IB(\ctc^{(0,0)}_{\infty,X}))\cdots\to
\Sect(U,\IB(\ctc^{(0,n)}_{\infty,X}))\to 0  
\end{equation}
is strictly exact in $\Ind(\ctBan)$.  For any open ball $b$ of $X$, 
Cartan's Theorem B shows that
$$
H^{k}(b,\O_{X})\simeq 0 \qquad (k>0).
$$
Hence, the sequence
$$
0\to\Sect(b,\IB(\O_{X}))\to
\Sect(b,\IB(\ctc^{(0,0)}_{\infty,X}))\cdots\to
\Sect(b,\IB(\ctc^{(0,n)}_{\infty,X}))\to 0  
$$
is strictly exact in $\Ind(\ctBan)$. Filtering inductive limits being 
exact in $\Ind(\ctBan)$, we see that
$$
0\to\IB(\O_{X})\to\IB(\ctc^{(0,0)}_{\infty,X})\to
\IB(\ctc^{(0,1)}_{\infty,X})\cdots\to
\IB(\ctc^{(0,n)}_{\infty,X})\to 0
$$
is a strictly exact sequence of sheaves with values in $\Ind(\ctBan)$.  
Moreover, since, by Proposition~\ref{prp:cinftyacyclic}, 
$\IB(\ctc^{(p,q)}_{\infty,U})$ is $\Sect(U,\cdot)$-acyclic, 
$\RSect(U,\IB(\O_{X}))$ is given by
$$
0\to\Sect(U,\IB(\ctc^{(0,0)}_{\infty,X}))\to
\Sect(U,\IB(\ctc^{(0,1)}_{\infty,X}))\cdots\to
\Sect(U,\IB(\ctc^{(0,n)}_{\infty,X}))\to 0.
$$
The sequence (*) being strictly exact, we get
$$
\RSect(U,\IB(\O_{X}))\simeq\Sect(U,\IB(\O_{X})).
$$
\end{proof}

\begin{proposition}\label{prp:sectacyclicibint}
If $X$ is a Stein manifold and $K$ is a holomorphically convex compact 
subset of $X$, we have
$$
\RSect(K,\IB(\O_{X}))\simeq\IB(\O_{X}(K)).
$$
\end{proposition}

\begin{proof}
It is well-known that $K$ has a fundamental system $\ctv$ of Stein 
open neighborhoods. By tautness, it follows that for $k>0$, we have
$$
LH^{k}(K,\IB(\O_{X}))\simeq
\indlim_{V\in\ctv}LH^{k}(V,\IB(\O_{X}))\simeq 0
$$
where the second isomorphism follows from 
Theorem~\ref{thm:ibsectacyclic}.  Hence, using 
Proposition~\ref{prp:ibisosect}, we get
$$
\RSect(K,\IB(\O_{X}))\simeq\Sect(K,\IB(\O_{X}))
\simeq\IB(\O_{X}(K)).
$$
\end{proof}

\begin{remark}\label{rem:omegaext}
Note that all the results in this section clearly hold if we replace 
$\O_{X}$ by the sheaf of holomorphic sections of holomorphic vector 
bundle. In particular, they hold for the sheaf $\Omega_{X}^{p}$ of 
holomorphic $p$-forms.
\end{remark}

\section{A factorization formula for $\IB(\O_{X\times Y})$}

\begin{definition}
For any $\rho=(\rho_{1},\cdots,\rho_{p})\in\ooint{0,+\infty}^{p}$, we 
set
$$
\Delta_{\rho}=\{z\in\C^{p}:|z_{1}|<\rho_{1},\cdots,|z_{p}|<\rho_{p}\}
$$
and we denote by $A_{\rho}$ the object of $\ctTc$ defined by endowing
$$
A_{\rho}=\{(a_{\alpha})_{\alpha\in\N^{p}}:
\sum_{\alpha}|a_{\alpha}|\rho^{\alpha}<+\infty\}
$$
with the norm
$$
\norm{(a_{\alpha})_{\alpha\in\N^{p}}}=
\sum_{\alpha}|a_{\alpha}|\rho^{\alpha}.
$$
\end{definition}

\begin{lemma}\label{lem:isoalone}
For any $\rho\in\ooint{0,+\infty}^{p}$, we have the isomorphism
$$
A_{\rho}\simeq l^{1}(\N^{p}).
$$
In particular, $A_{\rho}$ is a Banach space.
\end{lemma}

\begin{proof}
This follows directly from the fact that the application
$$
u:A_{\rho}\to l^{1}(\N^{p})
$$
defined by $u((a_{\alpha})_{\alpha\in\N^{p}})=
(a_{\alpha}\rho^{\alpha})_{\alpha\in\N^{p}}$ is continuous and 
bijective.
\end{proof}

\begin{lemma}\label{lem:isoadelta}
For any $p\in\N$,
$$
\indlim_{\rho\in\ooint{0,+\infty}^{p}}\IB(\O_{\C^{p}}(\Delta_{\rho}))
\simeq\indlim_{\rho\in\ooint{0,+\infty}^{p}}\IB(A_{\rho}).
$$
\end{lemma}

\begin{proof}
This follows directly from the fact that the canonical restriction 
morphism
$$
\O_{\C^{p}}(\Delta_{\rho'})\to\O_{\C^{p}}(\Delta_{\rho})
$$
may be factored through $A_{\rho}$ for $\rho'>\rho$.
\end{proof}

\begin{proposition}\label{prp:etenstopibplus}
Assume $X$, $Y$ are complex analytic manifolds. Then, there is a 
canonical isomorphism
$$
\IB(\O_{X})\letenstop\IB(\O_{Y})\simeq\IB(\O_{X\times Y}).
$$
\end{proposition}

\begin{proof}
Let $U$, $V$ be open subsets of $X$ and $Y$. The map
$$
u_{U,V}:\O_{X}(U)\times\O_{Y}(V)\to\O_{X\times Y}(U\times V)
$$
defined by setting
$$
u_{U,V}(f,g)(u,v)=f(u)g(v)
$$
is clearly bilinear and continuous. Hence, it induces a morphism
$$
\O_{X}(U)\tens_{\pi}\O_{Y}(V)\to\O_{X\times Y}(U\times V)
$$
and by Proposition~\ref{prp:bilintenstpcom}, we get a morphism
$$
\mu_{U,V}:\IB(\O_{X}(U))\tenstop\IB(\O_{Y}(V))\to
\IB(\O_{X\times Y}(U\times V))
$$
which is clearly well-behaved with respect to the restriction of 
$U$ or $V$. Therefore, we get a canonical morphism
$$
\mu:\IB(\O_{X})\etenstop\IB(\O_{Y})\to\IB(\O_{X\times Y}).
$$
To show that it is an isomorphism, it is sufficient to work at the 
level of germs and to prove that
$$
\mu_{(x,y)}:\IB(\O_{X})_{x}\tenstop\IB(\O_{Y})_{y}
\to\IB(\O_{X\times Y})_{(x,y)}
$$
is an isomorphism. The problem being local, we may assume $X=\C^{p}$, 
$Y=\C^{p'}$, $x=0$, $y=0$. In this case, Lemma~\ref{lem:isoadelta} 
shows that
$$
\IB(\O_{X})_{x}\simeq
\indlim_{\rho\in\ooint{0,+\infty}^{p}}\IB(A_{\rho}),\qquad
\IB(\O_{Y})_{y}\simeq
\indlim_{\rho'\in\ooint{0,+\infty}^{p'}}\IB(A_{\rho'})
$$
and
$$
\IB(\O_{X\times Y})_{(x,y)}\simeq
\indlim_{(\rho,\rho')\in\ooint{0,+\infty}^{p+p'}}
\IB(A_{(\rho,\rho')}).
$$
A direct computation shows that through these isomorphisms 
$\mu_{x,y}$ corresponds to the inductive limit of the maps
$$
\tau_{\rho,\rho'}:\IB(A_{\rho})\tenstop\IB(A_{\rho'})
\to\IB(A_{(\rho,\rho')})
$$
associated to the continuous bilinear maps
$$
t_{\rho,\rho'}:A_{\rho}\times A_{\rho'}\to A_{(\rho,\rho')}
$$
defined by
$$
t_{\rho,\rho'}((a_{\alpha})_{\alpha\in\N^{p}},
(a'_{\alpha'})_{\alpha'\in\N^{p'}})=
(a_{\alpha}a'_{\alpha'})_{(\alpha,\alpha')\in\N^{p+p'}}.
$$
Since the diagram
$$
\xymatrix{
\IB(A_{\rho})\tenstop\IB(A_{\rho'})
   \isoar[r] \ar[d]_{\tau_{\rho,\rho'}}
    &\cproobj{A_{\rho}\tenstop A_{\rho'}}
   \ar[d]^{\cproobj{t_{\rho,\rho'}}}
\\
\IB(A_{(\rho,\rho')}) \isoar[r]
    &\cproobj{A_{(\rho,\rho')}}
}
$$
is clearly commutative, to prove that $\mu_{(x,y)}$ is an isomorphism, 
it is sufficient to prove that $t_{\rho,\rho'}$ is an isomorphism.  
Thanks to Lemma~\ref{lem:isoalone}, this fact is an easy consequence 
of the well-known isomorphism
$$
l^{1}(\N^{p})\tenstop l^{1}(\N^{p'})\simeq 
l^{1}(\N^{p+p'}).
$$
By Proposition~\ref{prp:ibisosect},
$$
\IB(\O_{X})_{x}\simeq
\Sect(\{x\},\IB(\O_{X}))\simeq
\IB(\O_{X}(\{x\})).
$$
Since $\O_{X}(\{x\})$ is a DFN space, 
Proposition~\ref{prp:dfntensacyclic}, shows that
$$
\IB(\O_{X})_{x}\ltenstop\IB(\O_{Y})_{y}\simeq
\IB(\O_{X})_{x}\tenstop\IB(\O_{Y})_{y}.
$$
Therefore,
$$
\IB(\O_{X})\letenstop\IB(\O_{Y})\simeq
\IB(\O_{X})\etenstop\IB(\O_{Y})\simeq
\IB(\O_{X\times Y})
$$
as requested.
\end{proof}

\begin{corollary}\label{cor:holkun}
If $A$, $B$ are subsets of $X$ and $Y$ then
$$
\RSect_{c}(A\times B,\IB(\O_{X\times Y}))\simeq
\RSect_{c}(A,\IB(\O_{X}))\ltenstop\RSect_{c}(B;\IB(\O_{Y})).
$$
In particular, if $X$, $Y$ are Stein manifolds and $K$, $L$ are 
holomorphically convex compact subsets of $X$ and $Y$, then
$$
\IB(\O_{X\times Y}(K\times L))\simeq
\IB(\O_{X}(K))\tenstop\IB(\O_{Y}(L)).
$$
\end{corollary}

\begin{proof}
The first part is a direct consequence of 
Theorem~\ref{prp:etenstopibplus} and the K\"{u}nneth theorem for 
sheaves with values in $\Ind(\ctBan)$.  The second part follows from 
the first using Proposition~\ref{prp:sectacyclicibint}, 
Proposition~\ref{prp:dfntensacyclic} and the fact that $\O_{X}(K)$ is 
a DFN space.
\end{proof}

\section{Poincar\'{e} duality for $\IB(\O_{X})$}

\begin{proposition}
Assume $X$, $Y$ are complex analytic manifolds of dimension $d_{X}$ 
and $d_{Y}$.  Then, there is a canonical integration morphism
$$
\int_{X}:\reim{q_{Y}}(\IB(\Omega_{X\times Y})[d_{X\times Y}])
\to\IB(\Omega_{Y})[d_{Y}].
$$
\end{proposition}

\begin{proof}
Recall that integration along the fibers of $q_{Y}$ (i.e.~on $X$) 
defines morphisms
\begin{equation}\tag{*}
\int_{X}:\eim{q_{Y}}(\ctc^{p+d_{X},q+d_{X}}_{\infty,X\times Y})
\to\ctc^{p,q}_{\infty,Y}
\qquad(p,q\in\Z)
\end{equation}
which are compatible with $\partial$ and $\overline{\partial}$. Fix 
$p$, $q\in\Z$. Let $K$ be a compact subset of $X$ and let $U$ be an 
open subset of $Y$. One checks easily that the morphism
$$
\int_{X}:
\Sect_{K\times U}(X\times U;\ctc^{p+d_{X},q+d_{X}}_{\infty,X\times Y})
\to\Sect(U;\ctc^{p,q}_{\infty,Y})
$$
is continuous for the canonical topologies. Applying $\IB$, we get a 
morphism
$$
\Sect_{K\times U}(X\times U;
  \IB(\ctc^{p+d_{X},q+d_{X}}_{\infty,X\times Y}))
\to
\Sect(U;\IB(\ctc^{p,q}_{\infty,Y})).
$$
Taking the inductive limit on $K$, we get a morphism
$$
\Sect(U;\eim{q_{Y}}(\IB(\ctc^{p+d_{X},q+d_{X}}_{\infty,X\times Y})))
\to\Sect(U;\IB(\ctc^{p,q}_{\infty,Y}))
$$
and hence a morphism
$$
\eim{q_{Y}}(\IB(\ctc^{p+d_{X},q+d_{X}}_{\infty,X\times Y}))
\to\IB(\ctc^{p,q}_{\infty,Y})
$$
of sheaves with values in $\Ind(\ctBan)$. Thanks to the compatibility 
of (*) with $\partial$ and $\overline{\partial}$, we also get a 
morphism of complexes
$$
\eim{q_{Y}}(\IB(\ctc^{d_{X\times Y},\cdot}_{\infty,})
[d_{X\times Y}])\to
\IB(\ctc^{d_{X},\cdot}_{\infty, Y})[d_{Y}].
$$
Using the properties of Dolbeault resolutions, we get the requested 
integration morphism
$$
\int_{X}:\reim{q_{Y}}(\IB(\Omega_{X\times Y})[d_{X\times Y}])
\to\IB(\Omega_{Y})[d_{Y}].
$$
\end{proof}

\begin{remark}\label{rem:fubinilinearity}
Assume $X$, $Y$, $Z$ are complex analytic manifolds.  Then, one checks 
easily that Fubini Theorem gives rise to the commutative diagram
$$
\xymatrix{
\reim{q_{Z}}(\reim{q_{Y\times Z}}(\IB(\Omega_{X\times Y\times Z})
[d_{X\times Y\times Z}]))
        \ar[r]^-{\int_{X}} \isoar[d]
  & \reim{q_{Z}}(\IB(\Omega_{Y\times Z})[d_{Y\times Z}])
        \ar[d]^-{\int_{Y}}
\\
\reim{q_{Z}}(\IB(\Omega_{X\times Y\times Z})[d_{X\times Y\times Z}])
        \ar[r]_-{\int_{X\times Y}}
  & \IB(\Omega_{Z})[d_{Z}]
}
$$
Moreover, using the linearity of the integral, one gets the 
commutative diagram
$$
\xymatrix{
\reim{q_{Y}}(\IB(\Omega_{X\times Y})[d_{X\times Y}])
\tenstop\IB(\O_{Y})
          \ar[r]^-{\int_{X}\tenstop\id} \ar[d]_-{\text{projection}}
  & \IB(\Omega_{Y})[d_{Y}]\tenstop\IB(\O_{Y})
          \ar[dd]^-{\text{cup-product}}
\\
\reim{q_{Y}}(\IB(\Omega_{X\times Y})[d_{X\times Y}]
\tenstop q_{Y}^{-1}\IB(\O_{Y}))
          \ar[d]_{\text{cup-product}}
  &
\\
\reim{q_{Y}}(\IB(\Omega_{X\times Y})[d_{X\times Y}])
          \ar[r]_{\int_{X}}
  & \IB(\Omega_{Y})[d_{Y}]
}
$$
\end{remark}

\begin{theorem}\label{thm:dualibisoib}
Assume $X$ is a complex analytic manifold of dimension $d_{X}$ and 
denote $a_{X}:X\to\{\mathrm{pt}\}$ the canonical map.  Then, the 
morphism
$$
\IB(\Omega^{d_{X}-p}_{X})[d_{X}]\to\Dual(\IB(\Omega^{p}_{X})).
$$
induced by adjunction from
$$
\int_{X}\comp\smile:\eim{a_{X}}(\IB(\Omega^{d_{X}-p}_{X}[d_{X}])
\tenstop \IB(\Omega^{p}_{X}))\to\IB(\C)
$$
is an isomorphism.
\end{theorem}

\begin{proof}
The problem being local, it is sufficient to treat the case $p=0$ and 
to show that the morphism
$$
\RSect(U;\IB(\Omega_{U})[d_{U}])\to
\RHomtop(\RSect_{c}(U;\IB(\O_{U})),\IB(\C))
$$
obtained by adjunction from
$$
\int_{X}\comp\smile:\RSect(U;\IB(\Omega_{U})[d_{U}])
\ltenstop\RSect_{c}(U;\IB(\O_{U}))
\to\IB(\C)
$$
is an isomorphism for any open interval $U$ of $\C_{d_{U}}$.  This 
follows directly from Proposition~\ref{prp:isorsectrhomtop2} below 
with $V$ reduced to a point.
\end{proof}

\begin{remark}
As we will show elsewhere, the preceding theorem may be used to 
simplify the topological duality theory for coherent analytic sheaves.
\end{remark}

\begin{proposition}\label{prp:isorsectrhomtop2}
Assume $U$ is an open interval of $\C^{d_{U}}$ and $V$ is an open 
interval of $\C^{d_{V}}$.  Then, the canonical morphism
$$
\varphi_{U,V}:\RSect(U\times V,\IB(\Omega_{U\times V})[d_{U\times V}])
\to
\RHomtop(\RSect_{c}(U;\IB(\O_{U})),\RSect(V;\IB(\Omega_{V})[d_{V}]))
$$
obtained by adjunction from 
$$
\int_{X}\comp\smile:
\RSect(U\times V,\IB(\Omega_{U\times V})[d_{U\times V}])
\ltenstop\RSect_{c}(U;\IB(\O_{U}))
\to\RSect(V;\IB(\Omega_{V})[d_{V}])
$$
is an isomorphism.
\end{proposition}

\begin{proof}
Let $W$ be an open interval of $\C^{d_{W}}$ and assume that 
$\varphi_{U,V\times W}$ and $\varphi_{V,W}$ are isomorphisms.  Then, 
we have successively
\begin{align}
&\RSect(U\times V\times W;\IB(\Omega_{U\times V\times W})[d_{U\times 
  V\times W}])
\notag\\
&\quad\simeq \RHomtop(\RSect_{c}(U;\IB(\O_{U})),
  \RSect(V\times W,\IB(\Omega_{V\times W})[d_{V\times W}]))
\tag{1}\\
&\quad\simeq \RHomtop(\RSect_{c}(U;\IB(\O_{U})),
  \RHomtop(\RSect_{c}(V;\IB(\O_{V})),\RSect(W;\IB(\Omega_{W})[d_{W}])))
\tag{2}\\
&\quad\simeq \RHomtop(\RSect_{c}(U;\IB(\O_{U}))\ltenstop
  \RSect_{c}(V;\IB(\O_{V})),\RSect(W;\IB(\Omega_{W})[d_{W}]))
\tag{3}\\
&\quad\simeq\RHomtop(\RSect_{c}(U\times V;\IB(\O_{U\times V})),
  \RSect(W;\IB(\Omega_{W})[d_{W}])),
\tag{4}
\end{align}
where (1) and (2) follow from our assumptions, (3) is obtained by 
adjunction and (4) comes from Corollary~\ref{cor:holkun}.  Using 
Remark~\ref{rem:fubinilinearity}, we check easily that the composition 
of the preceding isomorphisms is equal to $\varphi_{U\times 
V,W}$.  Hence, an induction on $d_{U}$ reduces the problem to the case 
where $d_{U}=1$.  This will be dealt with in 
Proposition~\ref{prp:isorsectrhomtop} below.
\end{proof}

\begin{proposition}\label{prp:isorsectrhomtop}
Assume $U$ is an open interval of $\C$ and $V$ is an open interval 
of $\C^{n}$. Then, the canonical morphism
$$
\RSect(U\times V,\IB(\Omega_{U\times V})[d_{U\times V}])
\to
\RHomtop(\RSect_{c}(U;\IB(\O_{U})),\RSect(V;\IB(\Omega_{V})[d_{V}]))
$$
is an isomorphism.
\end{proposition}

\begin{proof}
For $P=\overline{U}$, sheaf theory gives us the two distinguished 
triangles
{\small
$$
\RSect_{\partial P\times V}(\C\times V,\IB(\Omega_{\C\times V}))\to
\RSect_{P\times V}(\C\times V,\IB(\Omega_{\C\times V}))\to
\RSect(U\times V,\IB(\Omega_{U\times V}))\to[+1]
$$
}%
and
$$
\RSect_{c}(U,\IB(\O_{U}))\to
\RSect(P,\IB(\O_{\C}))\to
\RSect(\partial P,\IB(\O_{\C}))\to[+1]
$$
where $\partial P$ denotes the boundary of $P$.  If we apply the 
functor $\RHomtop(\cdot,\IB(\Omega_{V}(V)))$ to the last triangle, we 
obtain the morphism of distinguished triangles
$$
\xymatrix{
\RSect_{\partial P\times V}(\C\times V,\IB(\O_{\C\times V}))[1]
    \ar[r]^-{\alpha}\ar[d]
  & \RHomtop(\RSect(\partial P,\IB(\O_{\C})),\IB(\Omega_{V}(V)))
    \ar[d]
\\
\RSect_{P\times V}(\C\times V,\IB(\O_{\C\times V}))[1]
    \ar[r]^-{\beta}\ar[d]
  & \RHomtop(\RSect(P,\IB(\O_{\C})),\IB(\Omega_{V}(V)))
    \ar[d]
\\
\RSect(U\times V,\IB(\O_{U\times V}))[1]
    \ar[r]^-{\gamma}\ar[d]^-{+1}
  & \RHomtop(\RSect_{c}(U,\IB(\O_{\C})),\IB(\Omega_{V}(V)))
    \ar[d]^-{+1}
\\
  &
}
$$
where $\alpha$ and $\beta$ are isomorphisms of the type considered in 
Proposition~\ref{prp:isofiniteunion} below ($\partial P$ is a finite 
union of closed intervals of $\C$).  It follows that $\gamma$ is an 
isomorphism.
\end{proof}

\begin{proposition}\label{prp:isofiniteunion}
Assume $K$ is a finite union of closed intervals of $\C$ and $V$ is 
an open interval of $\C^{n}$. Then, the canonical morphism
$$
\RSect_{K\times V}(\C\times V;\IB(\Omega_{\C\times V})[d_{\C\times V}])
\to
\RHomtop(\RSect(K;\IB(\O_{\C})),\RSect(V;\IB(\Omega_{V})[d_{V}]))
$$
obtained by adjunction from
$$
\int_{\C}\comp\smile:
\RSect_{K\times V}(\C\times V;\IB(\Omega_{\C\times V})[d_{\C\times V}])
\ltenstop\RSect(K;\IB(\O_{\C}))
\to\RSect(V;\IB(\Omega_{V})[d_{V}])
$$
is an isomorphism.
\end{proposition}

\begin{proof}
Assume first that $K$ is a closed interval of $\C$.  Since $P\times V$ 
is closed in $\C\times V$, we have the distinguished triangle 
{\small
$$
\RSect_{P\times V}(\C\times V,\IB(\O_{\C\times V}))\to
\RSect(\C\times V,\IB(\O_{\C\times V}))\to
\RSect((\C\setminus P)\times V,\IB(\O_{\C\times V}))\to[+1]
$$
}%
By Cartan's Theorem B and Theorem~\ref{thm:ibsectacyclic}, we have the 
isomorphisms
$$
\RSect(\C\times V,\IB(\O_{\C\times V}))\simeq
\IB(\O_{\C\times V}(\C\times V))
$$
and
$$
\RSect((\C\setminus P)\times V,\IB(\O_{\C\times V}))\simeq
\IB(\O_{\C\times V}((\C\setminus P)\times V)).
$$
Hence, the long exact sequence associated to the preceding 
distinguished triangle ensures that
$$
LH^{k}_{P\times V}(\C\times V,\IB(\O_{\C\times V}))=0
\qquad\forall k\geq 2
$$
and that the sequence
$$
\xymatrix{
0 \ar[r]
  & LH^{0}_{P\times V}(\C\times V,\IB(\O_{\C\times V})) \ar[r]
  & \IB(\O_{\C\times V}(\C\times V))  
      \ar`[r]`[ll]`[dll]`[dl][dl]
  & 
\\
  & IB(\O_{\C\times V}((\C\setminus P)\times V)) \ar[r]
  & LH^{1}_{P\times V}(\C\times V,\IB(\O_{\C\times V})) \ar[r]
  & 0
}
$$
is strictly exact.  Applying the functor $\IB$ to the sequence of 
Proposition~\ref{prp:splittc} below, we get the split exact sequence 
{\small
$$
0\to\IB(\O_{\C\times V}(\C\times V))\to
\IB(\O_{\C\times V}((\C\setminus P)\times V))\to
\IB(\Homtop[b](\O_{\C}(P),\O_{V}(V)))\to 0
\quad (*)
$$
}%
in $\Ind(\ctBan)$. Therefore, 
$$
LH^{0}_{P\times V}(\C\times V,\IB(\O_{\C\times V}))=0
$$
and
$$
LH^{1}_{P\times V}(\C\times V,\IB(\O_{\C\times V}))
\simeq\IB(\Homtop[b](\O_{\C}(P),\O_{V}(V))).
$$
Combining these results with Proposition~\ref{prp:isoibhomtop}, 
Theorem~\ref{thm:dfnfrhomtopacyc}, Theorem~\ref{thm:ibsectacyclic} and 
Proposition~\ref{prp:sectacyclicibint}, we obtain successively
\begin{align*}
\RSect_{P\times V}(\C\times V,\IB(\O_{\C\times V}))
  &\simeq\Homtop(\IB(\O_{\C}(P)),\IB(\O_{V}(V)))[-1]
\\
  &\simeq\RHomtop(\RSect(P;\IB(\O_{\C})),
    \RSect(V;\IB(\O_{V})))[-1].
\end{align*}
Thanks to Proposition~\ref{prp:intkothe} below, it follows easily that 
the canonical morphism
$$
\RSect_{P\times V}(\C\times V;\IB(\Omega_{\C\times V})[d_{\C\times V}])
\to\RHomtop(\RSect(P;\IB(\O_{\C})),\RSect(V;\IB(\Omega_{V})[d_{V}]))
$$
is an isomorphism.

Assume now that the result has been established when $K$ is a union of 
$k<N$ closed intervals of $\C$ and let us prove it when
$$
K=\Union_{i=1}^{N}P_{i}
$$
where $P_{i}$ ($i=1,\cdots,N$) is a closed interval of $\C$.  Set 
$L=\Union_{i=1}^{N-1}P_{i}$ and $Q=P_{N}$.  By the Mayer-Vietoris 
theorem associated to the decomposition $K=L\union Q$, we have the 
distinguished triangle
$$
\xymatrix{
\RSect(K,\IB(\O_{\C}))
    \ar[d]
\\
\RSect(L,\IB(\O_{\C}))\dsum\RSect(Q,\IB(\O_{\C}))
    \ar[d]
\\
\RSect(L\inter Q,\IB(\O_{\C}))
    \ar[d]^{+1}
\\
\null
}
$$
Applying the functor $\RHomtop(\cdot,\IB(\Omega_{V}(V)))$, we obtain 
the distinguished triangle
$$
\xymatrix{
A=\RHomtop(\RSect(L\inter Q,\IB(\O_{\C})),\IB(\Omega_{V}(V)))
    \ar[d]
\\
B=\RHomtop(\RSect(L,\IB(\O_{\C}))\dsum\RSect(Q,
    \IB(\O_{\C})),\IB(\Omega_{V}(V)))
    \ar[d]
\\
C=\RHomtop(\RSect(K,\IB(\O_{\C})),\IB(\Omega_{V}(V)))
    \ar[d]^{+1}
\\
\null
}
$$
Now, consider the Mayer-Vietoris distinguished triangle
$$
\xymatrix{
A'=\RSect_{(L\inter Q)\times V}
    (\C\times V,\IB(\Omega_{\C\times V}))
    \ar[d]
\\
B'=\RSect_{L\times V}(\C\times V,\IB(\Omega_{\C\times V}))\dsum
    \RSect_{Q\times V}(\C\times V,\IB(\Omega_{\C\times V}))
    \ar[d]
\\
C'=\RSect_{K\times V}(\C\times V,\IB(\Omega_{\C\times V}))
    \ar[d]
\\
\null
}
$$
Since $L\inter Q=\Union_{i=1}^{N-1}(P_{i}\inter P_{N})$ is a union of 
$N-1$ closed intervals of $\C$, the canonical morphisms
$$
A'[1]\to A \qquad\text{and}\qquad
B'[1]\to B
$$
are isomorphisms. The canonical diagram
$$
\xymatrix{
A'[1]\ar[r]\isoar[d]
  & B'[1]\ar[r]\isoar[d]
  & C'[1]\ar[r]^{+1}\ar[d]
  &
\\
A\ar[r]
  & B\ar[r]
  & C\ar[r]^{+1}
  &
}
$$
being commutative, the canonical morphism
$$
C'[1]\to C
$$
is also an isomorphism and the conclusion follows.
\end{proof}

\begin{proposition}\label{prp:intkothe}
Let $P$ be a compact interval of $\C$ and let $V$ be an open interval 
of $\C^{n}$. Then,
$$
\int_{\C}:H^{1}_{P\times V}(\C\times V,\Omega_{\C\times V})
\to H^{0}(V,\Omega_{V})
$$
sends the class of
$$
\omega=h(z,v)dz\wedge dv\in H^{0}(\C\setminus P,\Omega_{\C\times V})
$$
to
$$
\left(\int_{\partial P'}h(z,v)dz\right)dv
$$
where $P'$ is a compact interval of $\C$ such that 
${P'}^{\circ}\supset P$.
\end{proposition}

\begin{proof}
Let $\cti^{\cdot}$ be an injective resolution of $\Omega_{\C\times 
V}$. Denote
$$
u^{\cdot}:\ctc^{(v+1,\cdot)}_{\infty,\C\times V}\to\cti^{\cdot}
$$
a morphism extending $\id:\Omega_{\C\times V}\to\Omega_{\C\times V}$. 
The class $c$ of $\omega$ in
$$
H^{1}_{P\times V}(\C\times V,\Omega_{\C\times V})\simeq
H^{1}(\Sect_{P}(\C\times V,\cti^{\cdot}))
$$
is represented by $d\sigma$ where $\sigma\in\Sect(\C\times 
V,\cti^{0})$ extends $u^{0}(\omega)\in\Sect((\C\setminus P)\times 
V,\cti^{0})$. Let $\varphi$ be a function of class $C_{\infty}$ on 
$\C$ equals to $1$ on $\overline{\C\setminus P'}$ and to $0$ on 
$P''$, $P'$ and $P''$ being compact intervals such that 
${P''}^{\circ}\supset P$, ${P'}^{\circ}\supset P''$. Then, it is 
clear that $u^{0}(\varphi\omega)\in\Sect(\C\times V,\cti^{0})$ and 
that
$$
\sigma-u^{0}(\varphi\omega)\in\Sect_{c\times V}(\C\times V,\cti^{0}).
$$
Therefore, $d\sigma$ and $du^{0}(\varphi)$ give the same class in 
$H^{1}(\Sect_{c\times V}(\C\times V,\cti^{\cdot}))$.  It follows that 
$c'$ corresponds to the class of $\overline{\partial}(\varphi\omega)$ 
in $H^{1}(\Sect_{c\times V}(\C\times 
V,\ctc^{(v+1,\cdot)}_{\infty,\C\times V}))$.  Since $c'$ represents 
the image of $c$ by the canonical map
$$
H^{1}_{P\times V}(\C\times V,\Omega_{\C\times V})\to
H^{1}_{c\times V}(\C\times V,\Omega_{\C\times V}),
$$
we see that
$$
\int_{\C}c=
\int_{\C}\overline{\partial}(\varphi\omega)=
\int_{\overline{P'\setminus P''}}\overline{\partial}\varphi\omega=
\int_{\partial P'}\varphi\omega-\int_{\partial P''}\varphi\omega=
\int_{\partial P'}\omega.
$$
Hence the conclusion.
\end{proof}

\begin{proposition}\label{prp:splittc}
Let $P$ be a closed interval of $\C$ and let $V$ be an open interval 
of $\C^{n}$.  Then, in $\ctTc$, we have a split exact sequence of the 
form
$$
0\to\O_{\C\times V}(\C\times V)\to[r]
\O_{\C\times V}((\C\setminus P)\times V)\to[T]
\Homtop[b](\O_{\C}(P),\O_{V}(V))\to 0.
$$
where $r$ is the canonical restriction map and $T$ is defined by 
setting
$$
T(h)(\varphi)(v)=\int_{\partial P'}h(z,v)g(z)dz
$$
where $g$ is a holomorphic extension of $\varphi\in\O_{\C}(P)$ on an 
open neighborhood $U$ of $P$ and $P'$ is a compact interval of $\C$ 
such that ${P'}^{\circ}\supset P$ and $P'\subset U$.
\end{proposition}

\begin{proof}
Note that the definition of $T$ is meaningful since the right hand 
side clearly does not depend on the choices of $U$, $g$ and $P'$.  It 
is also clear that the function $T(h)(\varphi)$ is holomorphic on $V$ 
and that the operator $T$ is linear.  Let us show that $T$ is 
continuous.  Let $p$ be a continuous semi-norm of 
$\Homtop[b](\O_{\C}(P),\O_{V}(V))$.  We may assume that there is a 
bounded subset $B$ of $\O_{\C}(P)$ and a compact subset $K$ of $V$ 
such that
$$
p(\tau)=\sup_{\varphi\in B}\sup_{v\in K}|\tau(\varphi)(v)|,
\qquad \tau\in\Homtop[b](\O_{\C}(P),\O_{V}(V)).
$$
For $n>0$, set $U_{n}=\{u\in\C:d(u,P)<1/n\}$.  By cofinality, we have
$$
\O_{\C}(P)\simeq\indlim_{n>0}\O_{\C}(U_{n}).
$$
Moreover, for any $n>0$, $\O_{\C}(U_{n})$ is a Fr\'{e}chet space and 
the restriction
$$
\O_{\C}(U_{n})\to\O_{\C}(U_{n+1})
$$
is injective.  Hence, by~\cite[Chap.~IV, \S~19, 5.(5) 
(p.~225)]{Kothe69}, there is $n\in\N$ and a bounded subset $B_{n}$ of 
$\O_{\C}(U_{n})$ such that $B\subset r_{U_{n}}(B_{n})$.  Choosing a 
compact interval $P'_{n}$ of $\C$ such that ${P'_{n}}^{\circ}\supset P$ 
and $P'_{n}\subset U_{n}$, we see that
$$
p(T(h))\leq \sup_{g\in B_{n}}\sup_{v\in K}
\left|\int_{\partial P'_{n}}h(z,v)g(z)dz\right|
$$
and we can find $C>0$ such that
$$
p(T(h))
\leq C\sup_{g\in B_{n}}\sup_{z\in \partial P'_{n}}|g(z)|
\sup_{(z,v)\in\partial P'_{n}\times K}|h(z,v)|.
$$

Let us consider the linear map
$$
S:\Homtop[b](\O_{\C}(P),\O_{V}(V))\to
\O_{\C\times V}((\C\setminus P)\times V)
$$
defined by setting
$$
S(\tau)(z,v)=\frac{1}{2i\pi}\tau\left(\frac{1}{z-u}\right)(v).
$$
Let us check that $S$ is continuous.  Consider a compact subset $K$ of 
$\C\setminus P$ and a compact subset $L$ of $V$.  The set
$$
B_{K}=\{\frac{1}{z-u}:z\in K\}
$$
being bounded in $\O_{\C}(\C\setminus K)$, 
$r_{P,\C\setminus K}(B_{K})$ is a bounded subset of $\O_{\C}(P)$ and 
we have
$$
\sup_{(z,v)\in K\times L}|S(\tau)(z,v)|\leq \frac{1}{2\pi}
\sup_{f\in r_{P,\C\setminus K}(B_{K})}\sup_{v\in L}|\tau(f)(v)|.
$$
For any $\tau\in\Homtop[b](\O_{\C}(P),\O_{V}(V))$ and 
$\varphi\in\O_{\C}(P)$, there is an open $U$ of $\C$, containing $P$ 
and $g\in\O_{\C}(U)$ such that $\varphi=r_{U}(g)$.  Let $K$ be a 
closed interval included in $U$ and such that $K^{\circ}\supset P$ and 
let $\ctc$ be the oriented boundary of $K$.  For any $v\in V$, 
we have using the continuity of $\tau$ and Cauchy representation 
formula
\begin{align*}
T(S(\tau))(\varphi)(v)
  &=\frac{1}{2i\pi}\int_{\ctc}\tau\left(\frac{1}{z-u}\right)
    (v)g(z)dz
\\
  &=\tau\left(\frac{1}{2i\pi}
    \int_{\ctc}\frac{g(z)}{z-u}dz\right)(v)
\\
  &=\tau(g)(v)
   =\tau(\varphi)(v).
\end{align*}
It follows that $T\comp S=\id$ or, in other words, that $S$ is a 
section of $T$. 

Let us consider the continuous linear map
$$
R:\O_{\C\times V}((\C\setminus P)\times V)\to
\O_{\C\times V}(\C\times V)
$$
defined as follows.  Let $h\in\O_{\C\times V}((\C\setminus P)\times V)$ 
and $z\in\C$.  Consider $R>0$ such that
$$
z\in P_{R}^{\circ}=\{z:d(z,P)<R\}.
$$
Then, for any $v\in V$, we set
$$
R(h)(z,v)=\frac{1}{2i\pi}\int_{\ctc_{R}}\frac{h(u,v)}{u-z}du
$$
where $\ctc_{R}$ is the oriented boundary of $P_{R}$.  Since, for any 
$f\in\O_{\C\times V}(\C\times V)$ and any $(z,v)\in\C\times V$, we 
have
$$
R(r(f))(z,v)=
\frac{1}{2i\pi}\int_{\ctc_{R}}\frac{r(f)(u,v)}{u-z}du=
\frac{1}{2i\pi}\int_{\ctc_{R}}\frac{f(u,v)}{u-z}du=
f(z,v),
$$
we see that $R\comp r=\id$. The map $R$ is thus a retraction of $r$.

Thanks to a well-known result of homological algebra, the proof will 
be complete if we show that
$$
r\comp R+S\comp T=\id.
$$
To this end, consider $h\in\O_{\C\times V}((\C\setminus P)\times V)$ 
and $(z,v)\in(\C\setminus P)\times V$.  Fix $R>0$ such that $z\in 
P_{R}^{\circ}$ and denote $\ctc_{R}$ the oriented boundary of $P_{R}$.  
Let $\ctc$ be the oriented boundary of a closed interval $K\subset 
P_{R}$ such that $z\not\in K$ and $K^{\circ}\supset P$.  Denoting 
$\Gamma_{R}$ the oriented boundary of $P_{R}\setminus K^{\circ}$ and 
using Cauchy integral formula, we get
\begin{align*}
(r\comp R+S\comp T)(h)(z,v)   
  &=\frac{1}{2i\pi}\int_{\ctc_{R}}\frac{h(\xi,v)}{\xi-z}d\xi
    +\frac{1}{2i\pi}T(h)\left(\frac{1}{z-u}\right)(v)
\\
  &=\frac{1}{2i\pi}\int_{\ctc_{R}}\frac{h(\xi,v)}{\xi-z}d\xi  
    +\frac{1}{2i\pi}\int_{\ctc}\frac{h(\xi,v)}{z-\xi}d\xi
\\
  &=\frac{1}{2i\pi}\int_{\Gamma_{R}}\frac{h(\xi,v)}{\xi-z}d\xi
  \\
  &=h(z,v).
\end{align*}
\end{proof}

\begin{remark}
The preceding result is a slightly more precise form of a special case 
of the K\"{o}the-Grothendieck duality theorem (see~\cite{Kothe53} 
and~\cite{Grothendieck53a,Grothendieck53b}).
\end{remark}

\section{A holomorphic Schwartz' kernel theorem}

\begin{definition}
Let $X$ and $Y$ be complex analytic manifolds. We define 
$\Omega_{X\times Y}^{(r,s)}$ to be the subsheaf of 
$\Omega_{X\times Y}^{r+s}$ whose sections are the holomorphic 
differential forms that are locally a finite sum of forms of the type
$$
\omega_{i,j}dx_{i_{1}}\wedge\cdots\wedge dx_{i_{r}}
\wedge dy_{i_{1}}\wedge\cdots\wedge dy_{i_{s}}
$$
where $x$ and $y$ are holomorphic local coordinate systems on $X$ and 
$Y$.
\end{definition}

\begin{remark}\label{rem:omegaxy}
Clearly, $\Sect(W;\Omega_{X\times Y}^{(r,s)})$ has a canonical 
structure of FN space for any open subset $W$ of $X\times Y$.  
Therefore, using Proposition~\ref{prp:ibsheaf}, we see that 
$\IB(\Omega_{X\times Y}^{(r,s)})$ is a sheaf with value in 
$\Ind(\ctBan)$.  Moreover, using Proposition~\ref{prp:etenstopibplus}, 
one can check easily that
$$
\IB(\Omega_{X\times Y}^{(r,s)})\simeq
\IB(\Omega_{X}^{r})\letenstop\IB(\Omega_{Y}^{s}).
$$
\end{remark}

\begin{theorem}\label{thm:orelduality}
Assume $X$, $Y$ are complex analytic manifolds of dimension $d_{X}$, 
$d_{Y}$.  Then, we have a canonical isomorphism
$$
\IB(\Omega_{X\times Y}^{(d_{X}-r,s)})[d_{X}]\simeq
\rhomtop(q_{X}^{-1}\IB(\Omega_{X}^{r}),\epb{q_{Y}}\IB(\Omega_{Y}^{s})).
$$
\end{theorem}

\begin{proof}
We have successively
\begin{align}
\rhomtop(q_{X}^{-1}\IB(\Omega_{X}^{r}),
         \epb{q_{Y}}\IB(\Omega_{Y}^{s})[d_{Y}])
  &\simeq\rhomtop(q_{X}^{-1}\IB(\Omega_{X}^{r}),
      \epb{q_{Y}}\Dual(\IB(\Omega_{Y}^{d_{Y}-s})))
\tag{1}\\
  &\simeq\rhomtop(q_{X}^{-1}\IB(\Omega_{X}^{r}),
      \Dual(q_{Y}^{-1}\IB(\Omega_{Y}^{d_{Y}-s})))
\notag\\
  &\simeq\rhomtop(q_{X}^{-1}\IB(\Omega_{X}^{r}),
    \rhomtop(q_{Y}^{-1}\IB(\Omega_{Y}^{d_{Y}-s}),\omega_{X\times Y}))
\notag\\
  &\simeq\rhomtop(\IB(\Omega_{X}^{r})
      \letenstop\IB(\Omega_{Y}^{d_{Y}-s}),\omega_{X\times Y})
\notag\\
  &\simeq\rhomtop(\IB(\Omega^{(r,d_{Y}-s)}_{X\times Y}),
                  \omega_{X\times Y})
\tag{2}\\
  &\simeq\Dual(\IB(\Omega^{(r,d_{Y}-s)}_{X\times Y}))
\notag\\
  &\simeq\IB(\Omega^{(d_{X}-r,s)}_{X\times Y})[d_{X\times Y}]
\tag{3}
\end{align}
where $\omega_{X\times Y}$ denotes the dualizing complex on $X\times 
Y$ for sheaves with values in $\Ind(\ctBan)$.  Note that (1) and (3) 
follow from Theorem~\ref{thm:dualibisoib} and that (2) comes from 
Remark~\ref{rem:omegaxy}.
\end{proof}

As a consequence, we may now give Proposition~\ref{prp:isofiniteunion} 
and Proposition~\ref{prp:isorsectrhomtop} their full generality.

\begin{corollary}\label{cor:rsectisorhomtop}
Let $X$, $Y$ be complex analytic manifolds of dimension $d_{X}$ and 
$d_{Y}$.  Assume $K$ is a compact subset of $X$.  Then,
$$
\RSect_{K\times Y}(X\times Y;
\IB(\Omega_{X\times Y}^{(d_{X}-r,s)})[d_{X}])
\simeq
\RHomtop(\RSect(K;\IB(\Omega_{X}^{r}));
\RSect(Y;\IB(\Omega_{Y}^{s}))).
$$
Moreover, if $X$ and $Y$ are Stein manifolds and $K$ is 
holomorphically convex in $X$, these complexes are concentrated in 
degree $0$ and isomorphic to
$$
\IB(\Homtop[b](\Omega_{X}^{r}(K),\Omega_{Y}^{s}(Y))).
$$
\end{corollary}

\begin{proof}
Transposing to sheaves with values in $\Ind(\ctBan)$ a classical 
result of the theory of abelian sheaves, we see that
$$
\RSect_{K\times Y}(X\times Y;
\rhomtop(\opb{q_{X}}\shf,\epb{q_{Y}}\shg))
\simeq
\RHomtop(\RSect(K;\shf);\RSect(Y;\shg))
$$
if $\shf$ and $\shg$ are objects of $\Shv(X;\Ind(\ctBan))$ and 
$\Shv(Y;\Ind(\ctBan))$.  This formula combined with 
Theorem~\ref{thm:orelduality} gives the first part of the result.  The 
second part follows from Proposition~\ref{prp:sectacyclicibint}, 
Theorem~\ref{thm:ibsectacyclic} (using Remark~\ref{rem:omegaext}, 
Theorem~\ref{thm:dfnfrhomtopacyc} and 
Proposition~\ref{prp:isoibhomtop}.
\end{proof}

\begin{corollary}
Let $X$, $Y$ be complex analytic manifolds of dimension $d_{X}$ and 
$d_{Y}$.  Then,
$$
\RSect(X\times Y;\IB(\Omega_{X\times Y}^{(d_{X}-r,s)})[d_{X}])\simeq
\RHomtop(\RSect_{c}(X;\IB(\Omega_{X}^{r})),
         \RSect(Y;\IB(\Omega_{Y}^{s}))).
$$
\end{corollary}

\begin{proof}
This follows directly from the general isomorphism
$$
\RSect(X\times Y;\rhomtop(\opb{q_{X}}\shf,\epb{q_{Y}}\shg))\simeq
\RHomtop(\RSect_{c}(X;\shf),\RSect(Y;\shg))
$$
which holds for any objects $\shf$ and $\shg$ of 
$\Shv(X;\Ind(\ctBan))$ and $\Shv(Y;\Ind(\ctBan))$.
\end{proof}

\begin{lemma}\label{lem:submanvan}
Let $X$ be a complex analytic manifold of dimension $d_{X}$ and let 
$Y$ be a complex analytic submanifold of $X$ of dimension $d_{Y}$. 
Then,
$$
LH^{k}(\RSect_{Y}(\IB(\O_{X})))\simeq 0
$$
for $k\neq d_{X}-d_{Y}$.
\end{lemma}

\begin{proof}
Since the problem is local, it is sufficient to show that
$$
LH^{k}(\RSect_{\{0\}\times V}(U\times V;\IB(\O_{U\times V})))
\simeq 0
$$
for $k\neq d_{X}-d_{Y}$ if $U$ and $V$ are Stein open neighborhoods 
of $0$ in $\C^{d_{X}-d_{Y}}$ and $\C^{d_{Y}}$. In this situation, 
$\{0\}$ is a holomorphically convex compact subset of $U$ and we get 
from Corollary~\ref{cor:rsectisorhomtop} that
$$
\RSect_{\{0\}\times V}(U\times V;\IB(\O_{U\times V})[d_{X}-d_{Y}])
\simeq 
\IB(\Homtop[b](\O_{U}(\{0\}),\O_{V}(V))).
$$
The conclusion follows directly.
\end{proof}

\begin{theorem}
For any morphism of complex analytic manifolds $f:X\to Y$, we have a 
canonical isomorphism
$$
\rhomtop(\opb{f}\IB(\O_{Y}),\IB(\O_{X}))\simeq
\opb{\delta_{f}}\RSect_{\Delta_{f}}
\IB(\Omega_{X\times Y}^{(0,d_{Y})})[d_{Y}]
$$
where $\Delta_{f}$ is the graph of $f$ in $X\times Y$ and 
$\delta_{f}:X\to X\times Y$ is the associated graph embedding.  In 
particular,
$$
LH^{k}(\rhomtop(\opb{f}\IB(\O_{Y}),\IB(\O_{X})))=0
$$
for $k\neq 0$ and
$$
\rhom(\opb{f}\IB(\O_{Y}),\IB(\O_{X}))\simeq
\D_{X\to Y}^{\infty}.
$$
\end{theorem}

\begin{proof}
Using Theorem~\ref{thm:orelduality}, we see that 
$$
\IB(\Omega_{X\times Y}^{(0,d_{Y})})[d_{Y}]\simeq
\rhomtop(\opb{q_{Y}}\IB(\O_{Y}),\epb{q_{X}}\IB(\O_{X})).
$$
Applying $\epb{\delta_{f}}$, we get successively
\begin{align*}
\epb{\delta_{f}}\IB(\Omega_{X\times Y}^{(0,d_{Y})})[d_{Y}]
  &\simeq\epb{\delta_{f}}\rhomtop(\opb{q_{Y}}\IB(\O_{Y}),
    \epb{q_{X}}\IB(\O_{X}))
\\
  &\simeq\rhomtop(\opb{\delta_{f}}\opb{q_{Y}}\IB(\O_{Y}),
    \epb{\delta_{f}}\epb{q_{X}}\IB(\O_{X}))
\\
  &\simeq\rhomtop(\opb{(q_{Y}\comp\delta_{f})}\IB(\O_{Y}),
    \epb{(q_{X}\comp \delta_{f})}\IB(\O_{X}))
\\
  &\simeq\rhomtop(\opb{f}\IB(\O_{Y}),\IB(\O_{X})).
\end{align*}
This gives the first part of the result.  To get the second one, it is 
sufficient to use Lemma~\ref{lem:submanvan}, if we remember that, 
following~\cite{SKK73}, we have
$$
\D_{X\to Y}^{\infty}\simeq 
\opb{\delta_{f}}\RSect_{\Delta_{f}}
\Omega_{X\times Y}^{(0,d_{Y})}[d_{Y}].
$$
\end{proof}

\begin{corollary}\label{cor:rldinfty}
For any complex analytic manifold $X$ of dimension $d_{X}$, we have a 
canonical isomorphism
$$
\rhomtop(\IB(\O_{X}),\IB(\O_{X}))\simeq
\opb{\delta}\RSect_{\Delta}
\IB(\Omega_{X\times X}^{(0,d_{X})})[d_{X}]
$$
where $\Delta$ is the diagonal of $X\times X$ and 
$\delta:X\to X\times X$ is the diagonal embedding.  In 
particular,
$$
LH^{k}(\rhomtop(\IB(\O_{X}),\IB(\O_{X})))=0
$$
for $k\neq 0$ and
$$
\rhom(\IB(\O_{X}),\IB(\O_{X}))\simeq
\D_{X}^{\infty}.
$$
\end{corollary}

\begin{remark}
Note that the fact that continuous endomorphisms of $\O_{X}$ may be 
identified with partial differential operators of infinite order was 
conjectured by Sato and proved in~\cite{Ishimura78}. The vanishing of 
the topological $\ext{k}$ ($k>0$) is, to our knowledge, entirely new.
\end{remark}

\section{Reconstruction theorem}

Let $\ctr$ be a ring on $X$ with values in $\Ind(\ctBan)$ (i.e.\ a 
ring of the closed category $\Shv(X;\Ind(\ctBan))$ 
(see~\cite{Schneiders99})).  Denote by $\ctMod(\ctr)$ the 
quasi-abelian category formed by $\ctr$-modules.

If $\ctm$, $\ctn$ are two $\ctr$-modules, one sees easily that 
$\homtop(\ctm,\ctn)$ is endowed with both a structure of right 
$\ctr$-module and a compatible structure of left $\ctr$-module. These 
structures give two maps
$$
\xymatrix{
\homtop(\ctm,\ctn) \ar@<1ex>[r] \ar@<-1ex>[r]
  & \homtop(\ctr,\homtop(\ctm,\ctn)).
}
$$
As usual, we denote their equalizer by $\homtop[\ctr](\ctm,\ctn)$. In 
this way, we get a functor
$$
\homtop[\ctr](\cdot,\cdot):
\ctMod(\ctr)^{\op}\times\ctMod(\ctr)\to
\Shv(X;\Ind(\ctBan))
$$
which is clearly continuous on each variable and in particular left 
exact.  Using the techniques of~\cite{Schneiders99}, one sees easily 
that $\ctMod(\ctr)$ has enough injective objects and working as 
in~\cite[Proposition~2.3.10]{Schneiders99}, one sees that the functor 
$\homtop[\ctr](\cdot,\cdot)$ has a right derived functor
$$
\rhomtop[\ctr](\cdot,\cdot):
D^{-}(\ctMod(\ctr))^{\op}\times D^{+}(\ctMod(\ctr))
\to D^{+}(\Shv(X;\Ind(\ctBan))).
$$

Now, let $\cte$ be a sheaf on $X$ with values in $\Ind(\ctBan)$ and 
let $\ctn$ be an $\ctr$-module. Since $\homtop(\cte,\ctn)$ is 
canonically endowed with a structure of $\ctr$-module, we get a 
functor
$$
\homtop(\cdot,\cdot):
\Shv(X;\Ind(\ctBan))^{\op}\times\ctMod(\ctr)\to\ctMod(\ctr).
$$
One checks directly that this functor may be derived on the right by 
resolving the first argument by a complex of 
$K^{-}(\Shv(X;\Ind(\ctBan)))$ with components of the type
$$
\DSum_{i\in I}(P_{i})_{U_{i}}
$$
(where $P_{i}$ is a projective object of $\Ind(\ctBan)$ and $U_{i}$ is 
an open subset of $X$) and the second argument by a complex of 
$K^{+}(\ctMod(\ctr))$ with flabby components.  This gives us a derived 
functor
$$
\rhomtop(\cdot,\cdot):
D^{-}(\Shv(X;\Ind(\ctBan)))^{\op}\times D^{+}(\ctMod(\ctr))
\to D^{+}(\ctMod(\ctr))
$$
which reduces to the usual $\rhomtop$ functor if we forget the 
$\ctr$-module structures.

Finally, recall that an object $\ctm$ of $D^{b}(\ctMod(\ctr))$ is 
perfect if there are integers $p\leq q$ such that for any $x\in X$ 
there is a neighborhood $U$ of $x$ with the property that $\ctm|_{U}$ 
is isomorphic to a complex of the type
$$
0\to\shp^{p}\to\cdots\to\shp^{q}\to 0
$$
where each $\shp^{k}$ is a direct summand of a free $\shr_{U}$-module 
of finite type.  We denote by $D^{b}_{pf}(\ctMod(\ctr))$ the 
triangulated subcategory of $D^{b}(\ctMod(\ctr))$ formed by perfect 
objects.

\begin{proposition}\label{prp:isobimod}
Let $\ctn$ be a sheaf on $X$ with values in $\Ind(\ctBan)$ such that
$$
LH^{k}(\rhomtop(\ctn,\ctn))=0 \qquad(k\neq 0)
$$
and let $\ctr$ be the ring $\homtop(\ctn,\ctn)$ of internal 
endomorphisms of $\ctn$.  Then, $\ctn$ is an $\ctr$-module and the 
functor
$$
\rhomtop[\ctr](\cdot,\ctn):D^{b}_{pf}(\ctMod(\ctr))\to
D^{b}(\Shv(X;\Ind(\ctBan)))
$$
is well-defined. Moreover, we have a canonical isomorphism
$$
\rhomtop(\rhomtop[\ctr](\ctm,\ctn),\ctn)\simeq\ctm
$$
in $D(\ctMod(\ctr))$ for any $\ctm\in D^{b}_{pf}(\ctMod(\ctr))$. In 
particular, $\rhomtop[\ctr](\cdot,\ctn)$ identifies 
$D^{b}_{pf}(\ctMod(\ctr))$ with a full triangulated subcategory of 
$D^{b}(\Shv(X;\Ind(\ctBan)))$.
\end{proposition}

\begin{proof}
For any $\ctm\in D^{b}_{pf}(\ctMod(\ctr))$, it is clear that 
$$
\rhomtop[\ctr](\ctm,\ctn)\in D^{b}(\Shv(X;\Ind(\ctBan)))
$$
since $\rhomtop[\ctr](\ctr,\ctn)\simeq\ctn$.  The canonical morphism
$$
\ctm\ltenstop\rhomtop[\ctr](\ctm,\ctn)\to\ctn
$$
induces by adjunction a morphism
$$
\ctm\to\rhomtop(\rhomtop[\ctr](\ctm,\ctn),\ctn).
$$
If $\ctm\simeq\ctr$, $\rhomtop[\ctr](\ctm,\ctn)\simeq\ctn$ and
$$
\rhomtop(\rhomtop[\ctr](\ctm,\ctn),\ctn)\simeq
\rhomtop(\ctn,\ctn)\simeq\homtop(\ctn,\ctn)\simeq\ctr
$$
and the preceding morphism is an isomorphism.  It follows that it is 
also an isomorphism for $\ctm\simeq\ctr^{k}$ and, hence, if $\ctm$ is 
a direct summand of a free $\shr$-module of finite type.  Thanks to 
the local structure of perfect complexes, the conclusion follows 
easily.
\end{proof}

Let us consider the two functors
\begin{align*}
I_{\ctv}:\ctv & \to\Ind(\ctBan)  \\
E             & \mapsto\findlim_{\substack{F\subset E\\ 
                \dim F<+\infty}}F
\end{align*}
and
\begin{align*}
L_{\ctv}:\Ind(\ctBan)    & \to\ctv   \\
\findlim_{i\in\cti}E_{i} & \mapsto \indlim_{i\in\cti}E_{i}
\end{align*}
where $\ctv$ denotes the category of $\C$-vector spaces. They are 
clearly linked by the adjunction formula
$$
\Hom(I_{\ctv}(E),F)\simeq\Hom(E,L_{\ctv}(F))
$$
and they are both exact. Moreover,
$$
L_{\ctv}\comp I_{\ctv}=\id.
$$
For any sheaf $E$ on $X$ with values in $\ctv$, we denote 
$\tilde{I}_{\ctv}(E)$ the sheaf associated to the presheaf
$$
U\mapsto I_{\ctv}(E(U)).
$$
Similarly, to any sheaf $F$ on $X$ with values in $\Ind(\ctBan)$, we 
denote $\tilde{L}_{\ctv}(F)$ the sheaf
$$
U\mapsto L_{\ctv}(F(U))
$$
Working at the level of fibers, one checks easily that
$$
\tilde{L}_{\ctv}\comp\tilde{I}_{\ctv}=\id.
$$

\begin{proposition}\label{prp:rhomtoptilde}
Let $\ctn$ be a sheaf on $X$ with values in $\Ind(\ctBan)$ such that
$$
LH^{k}(\rhom(\ctn,\ctn))=0 \qquad(k\neq 0)
$$
and let $\ctr_{\ctv}$ be the ring $\hom(\ctn,\ctn)$ of endomorphisms 
of $\ctn$.  Then, $\ctn$ is an 
$\tilde{I}_{\ctv}(\ctr_{\ctv})$-module and the functor
$$
\rhomtop[\tilde{I}_{\ctv}(\ctr_{\ctv})](\tilde{I}_{\ctv}(\cdot),\ctn):
D^{b}_{pf}(\ctMod(\ctr_{\ctv}))\to D^{b}(\Shv(X;\Ind(\ctBan)))
$$
is well-defined. Moreover, we have a canonical isomorphism
$$
\rhom(
  \rhomtop[\tilde{I}_{\ctv}(\ctr_{\ctv})](
    \tilde{I}_{\ctv}(\ctm),\ctn),
    \ctn
)\simeq\ctm
$$
in $D(\ctMod(\ctr_{\ctv}))$ for any $\ctm\in 
D^{b}_{pf}(\ctMod(\ctr_{\ctv}))$.  

In particular, 
$\rhomtop[\tilde{I}_{\ctv}(\ctr_{\ctv})](\tilde{I}_{\ctv}(\cdot),\ctn)$ 
identifies $D^{b}_{pf}(\ctMod(\ctr_{\ctv}))$ with a full 
triangulated subcategory of $D^{b}(\Shv(X;\Ind(\ctBan)))$.
\end{proposition}

\begin{proof}
Applying $\tilde{L}_{\ctv}$ to the morphism
$$
\tilde{I}_{\ctv}(\ctm)\to
\rhomtop(\rhomtop[\tilde{I}_{\ctv}(\ctr_{\ctv})]
(\tilde{I}_{\ctv}(\ctm),\ctn),\ctn)
$$
we get a canonical morphism
$$
\ctm\to\rhom(\rhomtop[\tilde{I}_{\ctv}(\ctr_{\ctv})]
(\tilde{I}_{\ctv}(\ctm),\ctn),\ctn)
$$
since $\tilde{L}_{\ctv}\comp\rhomtop\simeq\rhom$.  The conclusion 
follows by working as in the proof of Proposition~\ref{prp:isobimod}.
\end{proof}

\begin{theorem}
Assume $X$ is a complex analytic manifold of dimension $d_{X}$.  Then, 
the sheaf $\IB(\O_{X})$ is an 
$\tilde{I}_{\ctv}(\D_{X}^{\infty})$-module and the functor
$$
\rhomtop[\tilde{I}_{\ctv}(\D_{X}^{\infty})]
(\tilde{I}_{\ctv}(\cdot),\IB(\O_{X})):
D^{b}_{pf}(\ctMod(\D_{X}^{\infty}))
\to D^{b}(\Shv(X;\Ind(\ctBan)))
$$
is well-defined. Moreover, we have a canonical isomorphism
$$
\rhom(\rhomtop[\tilde{I}_{\ctv}(\D_{X}^{\infty})]
(\tilde{I}_{\ctv}(\ctm),\IB(\O_{X})),
\IB(\O_{X}))\simeq\ctm
$$
in $D(\ctMod(\D_{X}^{\infty}))$ for any $\ctm\in 
D^{b}_{pf}(\ctMod(\D_{X}^{\infty}))$.  

In particular, $\rhomtop[\tilde{I}_{\ctv}(\D_{X}^{\infty})] 
(\tilde{I}_{\ctv}(\cdot),\IB(\O_{X}))$ identifies 
$D^{b}_{pf}(\ctMod(\D_{X}^{\infty}))$ with a full triangulated 
subcategory of $D^{b}(\Shv(X;\Ind(\ctBan)))$.
\end{theorem}

\begin{proof}
Thanks to Corollary~\ref{cor:rldinfty}, this is an easy consequence of 
Proposition~\ref{prp:rhomtoptilde}.
\end{proof}

\providecommand{\bysame}{\leavevmode\hbox to3em{\hrulefill}\thinspace}

\end{document}